\documentclass[usenames,dvipsnames, 11pt]{article}
\usepackage[utf8]{inputenc}

\usepackage[margin=2.5cm,a4paper]{geometry} 
\usepackage{multirow,listings,setspace,gnuplottex,latexsym,keyval,ifthen,moreverb,lscape,forest}
\usepackage{pgf,tikz}
\usetikzlibrary{arrows}

\usepackage{amsmath, amssymb, amsthm, bm}
\usepackage{thm-restate}
\usepackage{hyperref}
\usepackage[capitalize]{cleveref}
\usepackage{color}
\usepackage{url}
\usepackage{enumitem}
\usepackage{graphicx}
\usepackage{appendix}
\usepackage{cite}
\usepackage{mathtools}

\newtheorem{thm}{Theorem}[section]
\crefname{thm}{Theorem}{Theorems}
\newtheorem{prop}[thm]{Proposition}
\crefname{prop}{Proposition}{Propositions}
\newtheorem{lem}[thm]{Lemma}
\crefname{lem}{Lemma}{Lemmas}
\newtheorem{cor}[thm]{Corollary}
\newtheorem{clm}[thm]{Claim}

\newtheorem{fact}[thm]{Fact}
\numberwithin{thm}{section}
\numberwithin{equation}{section}

\theoremstyle{definition}
\newtheorem{mydef}[thm]{Definition}

\crefname{figure}{Figure}{Figures}

\DeclarePairedDelimiter{\parens}{(}{)}
\DeclarePairedDelimiter{\set}{\{}{\}}

\DeclarePairedDelimiter{\floor}{\lfloor}{\rfloor}
\DeclarePairedDelimiter{\ceil}{\lceil}{\rceil}
\DeclarePairedDelimiter{\brackets}{[}{]}
\DeclarePairedDelimiter{\size}{|}{|}

\def\cL{{\mathcal L}}
\def\eps{{\varepsilon}}

\title{Ramsey goodness of $k$-uniform paths, or the lack thereof}

\author{Simona Boyadzhiyska\thanks{School of Mathematics, University of Birmingham, Edgbaston, Birmingham, B15 2TT, UK.
 The research leading to these results was supported by EPSRC, grant no. EP/V002279/1 (A.~Lo) and EP/V048287/1 (S.~Boyadzhiyska and A.~Lo).
There are no additional data beyond that contained within the main manuscript.}\;\thanks{E-mail: \textsf{s.s.boyadzhiyska@bham.ac.uk}}
	\and
	Allan Lo\footnotemark[1]\;\thanks{E-mail: \textsf{s.a.lo@bham.ac.uk}}
}

\begin{document}
\maketitle

\begin{abstract}
    Given a pair of~$k$-uniform hypergraphs~$(G,H)$, the \emph{Ramsey number} of~$(G,H)$, denoted by~$R(G,H)$, is the smallest integer~$n$ such that in every red/blue-colouring of the edges of~$K_n^{(k)}$ there exists a red copy of~$G$ or a blue copy of~$H$.
Burr showed that, for any pair of graphs~$(G,H)$, where $G$ is large and connected, $R(G,H)\geq(v(G)-1)(\chi(H)-1)+\sigma(H)$, where~$\sigma(H)$ stands for the minimum size of a colour class over all proper~$\chi(H)$-colourings of~$H$. We say that~$G$ is \emph{$H$-good} if $R(G,H)$ is equal to the general lower bound.  Burr showed that, for any graph~$H$, every sufficiently long path is~$H$-good.

Our goal is to explore the notion of Ramsey goodness in the setting of $k$-uniform hypergraphs.  We demonstrate that, in stark contrast to the graph case, $k$-uniform $\ell$-paths are not~$H$-good for a large class of $k$-graphs. 
On the other hand, we prove that long loose paths are always at least \emph{asymptotically} $H$-good for every $H$ and derive lower and upper bounds that are best possible in a certain sense. 

In the 3-uniform setting, we complement our negative result with a positive one, in which we determine the Ramsey number asymptotically for pairs containing a long tight path and a 3-graph $H$ when~$H$ belongs to a certain family of hypergraphs. This extends a result of Balogh, Clemen, Skokan, and Wagner for the Fano plane asymptotically to a much larger family of 3-graphs.
\end{abstract}

\section{Introduction}

A~$k$-uniform hypergraph~$H$, or a~$k$-graph for short, consists of a (finite) set~$V(H)$ of vertices and a set~$E(H)$ of 
$k$-element subsets of~$V(H)$, called (hyper)edges.
Given~$k$-graphs~$G$ and~$H$, the \emph{Ramsey number} of the pair~$(G,H)$, denoted by~$R(G,H)$, is the smallest integer~$n$ such that, in every red/blue-colouring of the edges of the complete~$k$-graph~$K^{(k)}_n$ on $n$ vertices, we can find a 
red copy of~$G$ or a blue copy of~$H$.
Ramsey's seminal result~\cite{ramsey1930formallogic} implies that~$R(G,H)$ is finite for any pair of~$k$-graphs~$G$ and~$H$. Since then, the study of Ramsey numbers has become a prominent area of research in combinatorics and has inspired the development of many powerful tools in the field (see for example~\cite{conlon_2015_survey,mubayi2020survey} and the references therein).

Even in the simplest setting, when the uniformity is two, Ramsey numbers are often  notoriously difficult to understand. 
The most well-studied case is when~$G = H = K_t$. 
It is known from the early work of Erd\H{o}s~\cite{erdos1947} and Erd\H{o}s and Szekeres~\cite{erdos_combinatorial_1935} that, up to lower order terms, ~$2^{t/2}\leq R(K_t,K_t)\leq 2^{2t}$ as~$t\to\infty$; these bounds remained essentially best possible for several decades, until very recently Campos, Griffiths,  Morris, and Sahasrabudhe~\cite{campos2023exponential} announced the first exponential improvement in the upper bound. 

Apart from demonstrating the difficulty of understanding Ramsey numbers, this example shows that Ramsey numbers can grow very quickly compared to~$v(G)$ and~$v(H)$. 
It is then natural to ask: how \emph{small} can Ramsey numbers be? A general lower bound in the graph case was shown by Burr~\cite{burr1981ramsey}. As usual, we denote the chromatic number of a graph $H$ by $\chi(H)$, and we write $\sigma(H)$ for the smallest possible size of a colour class in a proper colouring of~$H$ using~$\chi(H)$ colours. Additionally, throughout the paper, we always assume that~$G$ is connected. 
 Following a slightly weaker observation by Chv\'atal and Harary~\cite{chvatal1972generalized}, Burr~\cite{burr1981ramsey} showed that, for any $G$ and $H$ with $v(G)\geq \sigma(H)$, we have 
\begin{align}\label{eq:Ramsey-goodness}
    R(G,H) \geq (v(G)-1)(\chi(H)-1)+\sigma(H).
\end{align} 
Indeed, colour the complete graph of order~$(v(G)-1)(\chi(H)-1)+\sigma(H)-1$ so that the red edges form~$\chi(H)$ cliques, one of order~$\sigma(H)-1$ and the rest of  order~$v(G)-1$; it is not difficult to check that there is neither a red copy of~$G$ nor a blue copy of~$H$ in this colouring. 
Classic results of Bondy and Erd\H{o}s\cite{bondy1973ramsey} and Chv\'atal~\cite{Chvtal1977TreecompleteGR}, predating Burr's work, show that the bound in~\eqref{eq:Ramsey-goodness} is attained with equality when the pair consists of a tree or a long cycle and a complete graph. Motivated by these early results,
Burr~\cite{burr1981ramsey} and Burr and Erd\H{o}s~\cite{burr1983generalizations} investigated what other pairs have this property, introducing the notion of \emph{Ramsey goodness}. More precisely, a graph~$G$ is said to be \emph{$H$-good} if the lower bound in~\eqref{eq:Ramsey-goodness} is attained for the pair~$(G,H)$. In this paper, we are interested in the following result due to Burr.

\begin{thm}[Burr~\cite{burr1981ramsey}]\label{thm:burr}
    For every graph $H$, there exists an integer $n_0 = n_0(H)$ such that every path or cycle on at least $n_0$ vertices is $H$-good.
\end{thm}  

Since its introduction Ramsey goodness has received considerable attention (see~\cite[Section~2.5]{conlon_2015_survey} and the references therein for some history and results). Typically in this line of research~$H$ is thought of as a fixed graph and the task is to identify what properties make a (sufficiently large) graph~$H$-good. 
Several conjectures were made (for example, by Burr~\cite{burr1981ramsey} and Burr and Erd\H{o}s~\cite{burr1983generalizations}), suggesting that, for a fixed graph~$H$, every sufficiently sparse large graph~$G$ should be~$H$-good. In support of these conjectures,  Burr, Erd\H{o}s,  Faudree, Rousseau, and Schelp\cite{burr1985ramsey}  proved that (large) bounded degree graphs are $H$-good for every bipartite graph~$H$. 
However, the conjectures turned out to be false in general, as shown by Brandt~\cite{brandt1996expanding}, whose result essentially proves that graphs with good expansion properties are not $H$-good for any non-bipartite $H$.
On the other hand, it is known that there are some families of graphs such that every sufficiently large member is~$H$-good for every~$H$. 
One example is given by \cref{thm:burr} above.
More generally, Allen, Brightwell, and Skokan~\cite{allen2013ramsey} proved that, for every fixed~$H$, every large graph with bounded bandwidth is $H$-good.

We are interested in exploring the notion of Ramsey goodness for hypergraphs. As usual, a proper colouring of a hypergraph~$H$ is a colouring of the vertices of~$H$ such that no edge of~$H$ is monochromatic; $\chi(H)$ is then defined as the minimum number of colours in a proper colouring of~$H$, and~$\sigma(H)$ is the smallest possible size of a colour class in a proper colouring of~$H$ using~$\chi(H)$ colours. Further, we say that a hypergraph is \emph{connected} if it is not a disjoint union of two smaller hypergraphs. We again assume that $G$ is a connected $k$-graph whenever necessary.
It is not difficult to check that Burr's argument proving~\eqref{eq:Ramsey-goodness} readily generalises to pairs of $k$-graphs $(G,H)$ satisfying  $v(G)\geq \sigma(H)$, and 
we say that~$G$ is \emph{$H$-good} if equality holds in \eqref{eq:Ramsey-goodness}.

In particular, we seek to determine to what extent \cref{thm:burr} generalises to $k$-graphs. For this, we first need to define a suitable notion of a path. Let $k\geq 3$, $\ell\geq 1$, and $q\geq 1$ be integers and $n = q(k-\ell) + \ell$. The $n$-vertex $k$-uniform \emph{$\ell$-path} $P^{(k)}_{n,\ell}$ consists of $n$ vertices $v_1,\dots,v_n$ and hyperedges $e_1,\dots,e_{q}$, where $e_i = \set{v_{(i-1)(k-\ell)+1},\dots, v_{i(k-\ell)+\ell}}$; in other words, an $\ell$-path consists of a sequence of $n$ vertices in which each edge is a subsequence of length $k$ and consecutive edges intersect in precisely $\ell$ vertices. The \emph{length} of an $\ell$-path is the number of edges it contains. 
In the special case where $\ell = k-1$, the corresponding $\ell$-path $P^{(k)}_{n,\ell}$ is called a \emph{tight path}, while when $\ell = 1$, we refer to $P^{(k)}_{n,\ell}$ as a \emph{loose path}. Notice that, while a tight path exists for every $n\geq k$, this is not true for a general $\ell$-path. To be precise,  an $\ell$-path on $n$ vertices exists if and only if $n\equiv \ell \pmod{k-\ell}$. Thus, when we talk about a $k$-uniform $\ell$-path on $n$ vertices, we will implicitly assume that $n$ satisfies this condition.

Similarly, given integers $k\geq 3$, $\ell\geq 1$, and $q\geq 1$ and $n = q(k-\ell)$, the $n$-vertex $k$-uniform \emph{$\ell$-cycle} $C^{(k)}_{n,\ell}$ consists of $n$ vertices $v_1,\dots,v_n$ and hyperedges $e_1,\dots,e_{q}$, where $e_i = \set{v_{(i-1)(k-\ell)+1},\dots, v_{i(k-\ell)+\ell}}$, where for convenience we set $v_{n+i} = v_i$ for all $i\in [\ell]$. Again, notice that an $\ell$-cycle on $n$ vertices exists if and only if $n\equiv 0 \pmod{k-\ell}$, and we will always assume that $n$ is of the correct form when referring to an $\ell$-cycle on $n$ vertices. Again, the \emph{length} of an $\ell$-cycle is the number of edges it contains.

The study of Ramsey goodness in hypergraphs was first undertaken by Balogh, Clemen, Skokan, and Wagner~\cite{balogh2020ramsey} and was motivated by a question of Conlon.  Letting~$\cL_2$ denote the \emph{Fano plane}, that is, the unique 3-graph on seven vertices in which every pair of vertices is contained in a unique edge, Conlon asked what 3-graphs are~$\cL_2$-good. Balogh, Clemen, Skokan, and Wagner~\cite{balogh2020ramsey} made progress towards answering this question by showing that any sufficiently long tight path is~$\cL_2$-good. 
\begin{thm}[Balogh, Clemen, Skokan, and Wagner~\cite{balogh2020ramsey}]\label{thm:bcsw}
    There exists an integer $n_0$ such that, for all $n\geq n_0$, the tight path $P_{n,2}^{(3)}$ is $\cL_2$-good. 
\end{thm}

In light of \cref{thm:burr,thm:bcsw}, we explore several natural directions.  One question we address is  what 3-graphs long tight paths are good for. Second, the Fano plane is a 3-graph arising from a projective plane, so we also study $R(P^{(q+1)}_{n,q},\cL_q)$, where $\cL_q$ denotes a hypergraph corresponding to a projective plane of order $q$. We also move away from tight paths and consider other types of $\ell$-paths, with a particular focus on loose paths.

%%%%%%%%%%%%%%%%%%%%%%%%%%%%%%%%%%%%%%%%%%%
\subsection{Results}

We first begin with a general result about $\ell$-paths that shows that \cref{thm:burr} does not extend to $\ell$-paths in higher uniformities when $\ell\geq 2$ even asymptotically. More specifically, for every~$k\geq 3$, we find a large class of $k$-graphs $H$ such that  $R(P^{(k)}_{n,\ell}, H)$ considerably exceeds the lower bound.

\begin{prop}\label{prop:ell-paths}
    Let $k\geq 3$, $\ell\geq 2$, and $n\geq 1$ be integers with $n\equiv \ell \pmod{k-\ell}$.  Let $H$ be a $k$-graph such that, for every proper colouring $A_1\cup \dots\cup A_{\chi(H)}$ of $H$ with colours in $[\chi(H)]$ and every colour $i\in [\chi(H)]$, there is an edge $e$ such that $\size{e\cap A_i} \geq 1$ and $\size{e\cap A_j} \leq \ell-1$ for all~$j\in [\chi(H)]\setminus\set{i}$.
    Then $R(P^{(k)}_{n,\ell}, H) \geq (\chi -1)(n-1) +\floor*{\frac{n}{k}}$.
\end{prop}

\subsubsection{Loose paths and cycles}

Considering \cref{prop:ell-paths}, it is natural to ask what happens in the remaining case, $\ell=1$. Our first result concerning loose paths shows that they are at least \emph{asymptotically} $H$-good for every $H$.

\begin{prop}\label{prop:loosepaths_asymptotic}
    For every $\eps>0$, every integer $k\geq 3$, and every $k$-graph $H$, there exists an integer $n_0 = n_0(k,H,\eps)$ such that $R(P^{(k)}_{n,1}, H) \leq (\chi(H)-1)n +\eps n$ for all $n\geq n_0$.
\end{prop}

As we will see below, however, loose paths are not always \emph{exactly} $H$-good. Nevertheless, we are able to show an upper bound that is very close to tight, and is best possible in a certain sense. Before we present these results, we define the following function that will appear in both our lower and upper bounds.  

\begin{mydef}\label{def:tau}
Given integers $k,\alpha \geq 1$, we define $\tau(k, \alpha)$ to be the largest integer $n$ such that there exists a $k$-graph on $n$ vertices with independence number less than $\alpha$ and no loose path of length two. 
\end{mydef}
Observe that, if $\alpha < k$, then $\tau(k, \alpha) = \alpha - 1$. In the other regime, we determine $\tau(k,\alpha)$ up to an additive constant of at most $k-2$.

\begin{prop} \label{prop:tau_bounds}
For integers $\alpha \ge k \ge 2$, we have $\alpha-1 +  (k-1)\floor*{\frac{\alpha-1}{k-1}}  \le  \tau (k, \alpha) \leq 2\alpha - 2$.
\end{prop}

We are now ready to state our general upper bound on the Ramsey number of a long loose path versus any $k$-graph $H$.

\begin{thm}\label{thm:loosepaths_upper}
Let $k \ge 3$ and $H$ be a $k$-graph. 
Then there exists an integer $n_0 = n_0(H)$ such that, for all $n \ge n_0$ satisfying $n \equiv 1 \pmod{k-1}$, we have 
\begin{align*}
R(P^{(k)}_{n,1}, H) \leq (\chi(H)-1) (n-1) + \max \{ \tau (k-1, \sigma(H))-2k+3, \sigma(H)\} .    
\end{align*} 
\end{thm}

Note that, when $\sigma(H)\leq 2k-1$, the loose path $P^{(k)}_{n,1}$ is $H$-good.  \cref{thm:loosepaths_upper} determines $R(P^{(k)}_{n,1}, H)$ up to an additive constant of \mbox{$\tau (k-1, \sigma(H))-2k+3$}, which by \cref{prop:tau_bounds} is at most $2\sigma(H)-2k+1$. Next we show that the upper bound in \cref{thm:loosepaths_upper} is best possible.

\begin{prop}\label{prop:loose_path_lower}
Let $k\geq 3$, $\chi \ge 2$, $t \geq k-1$, and $n \ge 3(k-1)$ be integers, and assume that $n \equiv 1 \pmod{k-1}$. 
Let $H$ be a $k$-graph on vertex set $V(H) = A_1\sqcup\dots \sqcup A_{\chi}$, where $\size{A_i} > (\chi-1)(k-2) + \tau(k-1,t)$ for all $i\in [\chi-1]$ and $\size{A_{\chi}} = t$, whose edges are all $k$-subsets of~$V(H)$ satisfying $\size{e \cap A_i} = k-1$ for some $i \in [\chi]$. 
Then $R(P^{(k)}_{n,1},H) \geq (\chi-1)(n-1) + \tau(k-1,t) -2k+3$.
\end{prop}

\medskip
Our methods allow us to derive a corresponding result for loose cycles.

\begin{thm}\label{thm:loosecycles_upper}
Let $k \ge 3$ and $H$ be a $k$-graph. 
Then there exists an integer $n_0 = n_0(H)$ such that, for all $n \ge n_0$ satisfying $n \equiv 
0\pmod{k-1}$, we have 
\begin{align*}
R(C^{(k)}_{n,1}, H) \leq (\chi(H)-1) (n-1) + 4k\binom{\chi(H)-1}{2} + \tau (k-1, \sigma(H)) +1.    
\end{align*} 
\end{thm}

Note that the additive constant  depends on $\chi(H)$. In the next result, we show that this dependence is necessary. We make no effort to optimise the constant.

\begin{prop}\label{prop:loose_cycle_lower}
Let $k\geq 3$, $t, q \geq k-1$, $\chi > \binom{q}{k-1}$, and $n \ge 3(k-1)$ be integers with $n \equiv 0 \pmod{k-1}$. 
Let $H$ be a $k$-graph on vertex set $V(H) = A_1\sqcup\dots \sqcup A_{\chi}$, where $\size{A_i} > (\chi-1)(k-2) + \max\set{\tau(k-1,t),q}$ for all $i\in [\chi-1]$ and $\size{A_{\chi}} = t$,  whose edges are all $k$-subsets of~$V(H)$ satisfying $\size{e \cap A_i} = k-1$ for some $i \in [\chi]$.
Then 
\begin{align*}
  R(C^{(k)}_{n,1},H) > (\chi-1)(n-1) + \max\set{\tau(k-1,t),q}.  
\end{align*}
\end{prop}

\subsubsection{Tight paths and cycles}

 A natural problem arising from the work of Balogh, Clemen, Skokan, and Wagner~\cite{balogh2020ramsey} (\cref{thm:bcsw}) is to study what happens in higher uniformities when the Fano plane is replaced by a higher-order projective plane. A \emph{projective plane~$\mathcal{L}_q$ of order~$q$} is a $(q+1)$-regular $(q+1)$-graph on $q^2+q+1$ vertices in which every pair of vertices is contained in a unique edge. Note that $\cL_2$ is the Fano plane. For $q\geq 3$, a result of Richardson~\cite{richardson1956finite} implies that $\chi(\cL_q) = 2$. \cref{prop:ell-paths} then implies that tight paths are not $\cL_q$-good.
\begin{cor}\label{prop:proj_planes}
    Let $\cL_q$ be any projective plane of order $q\geq 3$. Then, for any large integer $n$, we have $R(P^{(q+1)}_{n,q}, \cL_q) \geq \floor*{\frac{q+2}{q+1}(n-1)}$.
\end{cor}

The Fano plane is also a Steiner triple system, that is, a 3-graph in which every pair of vertices is contained in a unique edge. Thus, it is natural to ask whether long 3-uniform tight paths are $\mathcal{S}$-good for every Steiner triple system $\mathcal{S}$. Unfortunately, Forbes~\cite[Theorem~1.1]{forbes2003uniquely} showed the existence of infinitely many Steiner triple systems $\mathcal{S}$ with chromatic number 3 whose unique proper 3-colouring is equitable. Since every pair of vertices is contained in an edge in $\mathcal{S}$, \cref{prop:ell-paths} implies that $P_{n,2}^{(3)}$ is not $\mathcal{S}$-good.

Given \cref{prop:ell-paths}, we concentrate on 3-graphs~$H$ that have at least one~$\chi(H)$-colouring in which every edge intersects precisely two different colour classes. In fact, we restrict our attention to the following class of 3-graphs. For an integer~$\ell\geq 1$, we write~$T_{\ell}$ to denote a tournament on vertex set~$[\ell]$, that is, an orientation of the complete graph~$K_{\ell}$;~$TT_{\ell}$ stands for a transitive (i.e., acyclic) tournament on~$[\ell]$.

\begin{mydef}\label{def:tournament_hypergraph}
    Let~$\chi\geq 1$ be an integer and~$T_{\chi}$ be a tournament on~$[\chi]$. 
    We say that a 3-graph~$H$ is a \emph{tournament hypergraph associated to~$T_\chi$} if~$V(H)$ can be partitioned into sets~$A_1\cup\dots\cup A_{\chi}$ so that~$E(H) = \set{xyz: x,y\in A_i, z\in A_j, (i,j)\in E(T_\chi)}$, that is, the edge set of~$H$ consists of precisely those triples containing two vertices from some set~$A_i$ and a third vertex from some set~$A_j$, where~$(i,j)$ is an arc of~$T_{\chi}$. 
    For an integer~$m\geq 1$, we write~$H(T_{\chi},m)$ to denote a tournament hypergraph associated to~$T_{\chi}$ in which each vertex class~$A_i$ has size~$m$.
\end{mydef}

Let~$T_{\chi}$ be a non-transitive tournament and~$H = H(T_{\chi}, m)$. 
Surprisingly, in this case, not only are tight paths not $H$-good, but in fact~$R(P^{(3)}_{n,2},H)/n$
cannot be bounded above by any function depending only on~$\chi$. More precisely, the next proposition shows that $R(P^{(3)}_{n,2},H)\geq \frac{2}{3}(n-2)(m-1)$.

\begin{prop}\label{prop:non_transitive}
    Let~$\chi\geq 3$ and~$m\geq 2$ be integers,~$T_{\chi}$ be any non-transitive tournament, and~$H = H(T_{\chi},m)$. Let~$n,t\geq 1$ be integers such that~$\floor*{\frac{3t}{2}}+1 < n$. Then~$R(P^{(3)}_{n,2},H) \geq (m-1)t+1$.
\end{prop}

The situation is fairly different when~$H$ is a tournament hypergraph associated to a \emph{transitive} tournament~$TT_{\chi}$.  
Once again, the pair~$(P^{(3)}_{n,2},H)$ is generally not Ramsey good, but as we will see soon, in this case~$R(P^{(3)}_{n,2},H)$ \emph{can} be bounded above by a function depending only on~$\chi$ and~$n$. Given an integer~$\ell\geq 1$, let~$\vec{R}(\ell)$ denote the smallest integer~$N$ such that any tournament on at least~$N$ vertices contains a copy of~$TT_\ell$. 
It is well known that~$\vec{R}(\ell)$ is finite for any~$\ell\geq 1$.

\begin{thm}\label{thm:tight_paths_ub}
    Given integers~$\chi\geq 2$ and~$m\geq 2$ and a real number~$\eps>0$, there exists an integer~$n_0 =n_0(\chi,m,\eps)$ such that
    \[R(P^{(3)}_{n,2},H(TT_\chi,m)) \leq  \begin{cases}
    \parens*{1+\eps}n & \text{ if } \chi=2,\\
    \parens*{\frac{2}{3}+\eps}(\vec{R}(\chi)-1)n & \text{ if } \chi\geq 3,
    \end{cases}
    \]
    for all~$n \ge n_0$. 
\end{thm}

Notice that combining the well-known fact that~$\vec{R}(3) = 4$ with~\cref{thm:tight_paths_ub} shows that, if~$H$ has chromatic number three and is a subhypergraph of~$H(TT_3,m)$ for some integer~$m$, then~$R(P^{(3)}_{n,2},H)\leq \parens*{2+o(1)}n$ as~$n\to\infty$. Together with the lower bound from~\eqref{eq:Ramsey-goodness}, this result then determines~$R(P^{(3)}_{n,2},H)$ asymptotically for large~$n$. It is not difficult to check that the Fano plane satisfies these properties. Therefore the special case~$\chi=3$ of~\cref{thm:tight_paths_ub} proves  the result of Balogh, Clemen, Skokan, and Wagner~\cite{balogh2020ramsey} asymptotically, extending it to a large family of 3-graphs. In particular, if~$H$ belongs to this family, then~$P^{(3)}_{n,2}$ is \emph{asymptotically}~$H$-good as~$n\to\infty$.

It is known that $\vec{R}(\chi)$ grows exponentially with $\chi$, so the upper bound in~\cref{thm:tight_paths_ub} is considerably larger than the lower bound from~\eqref{eq:Ramsey-goodness}. Nevertheless, this upper bound is asymptotically tight.

\begin{prop}\label{prop:transitive_lb}
    Let~$\chi\geq 3$ and $m\geq\vec{R}(\chi)$ be integers.  
    Then~$H = H(TT_{\chi},m)$ satisfies~$R(P^{(3)}_{n,2}, H)\geq \parens*{\frac{2}{3}n-3}(\vec{R}(\chi)-1) + 1 = (1+o(1))\frac{2}{3}(\vec{R}(\chi)-1)n$ as~$n\to\infty$. 
\end{prop}

With minor modifications to the proofs we can prove an analogous result for tight cycles instead of tight paths. 

\begin{thm}\label{thm:tight_cycles_ub}
    Given integers~$\chi\geq 2$ and~$m\geq 2$ and a real number~$\eps>0$, there exists an integer~$n_0 =n_0(\chi,m,\eps)$ such that
    \[\parens*{\frac23n-3}(\vec{R}(\chi)-1) < R(C^{(3)}_{n,2},H(TT_\chi,m)) \leq  \begin{cases}
    \parens*{1+\eps}n & \text{ if } \chi=2,\\
    \parens*{\frac{2}{3}+\eps}(\vec{R}(\chi)-1)n & \text{ if } \chi\geq 3,
    \end{cases}
    \]
    for all~$n \ge n_0$.
\end{thm}

\subsection{Notation}
Our graph and hypergraph notation is mostly standard. 
Wherever necessary, we use superscripts to indicate the uniformity of a hypergraph. For a hypergraph $F$ and a subset $A\subseteq V(F)$, we write $F[A]$ for the subhypergraph of $F$ induced by $A$. A subset $I$ of vertices in a hypergraph $F$ is independent if $F[I]$ contains no edges; the independence number of $F$, denoted $\alpha(F)$, is the maximum size of an independent set in $F$. Given hypergraphs $H_1,H_2$, and $F$, a red/blue-edge-colouring of $F$ is $(H_1,H_2)$-free if there is no red copy of $H_1$ and no blue copy of $H_2$.

Let~$F$ be any 3-uniform hypergraph. Given three subsets~$A,B, C\subseteq V(F)$, we write~$E(A,B,C)$ for the set of ordered triples~$(a,b,c)\in A\times B\times C$ such that ~$a,b,$ and~$c$ are all distinct and~$abc\in E(F)$ (note that we do not require the sets~$A,B$, and~$C$ to be distinct or disjoint); for convenience, we sometimes refer to the edges in $E(A,B,C)$ as $ABC$-edges. We write~$e(A,B,C) = \size{E(A,B,C)}$ and~$d(A,B,C) = \frac{e(A,B,C)}{\size{\set{(a,b,c)\in (A,B,C): a\neq b,a\neq c, b\neq c}}}$. 
Further, if the edges of~$F$ are coloured red and blue, we write~$F_r$ and~$F_b$ for the red and blue subgraph of~$F$, respectively; we then define~$E_r(A,B,C), E_b(A,B,C), e_r(A,B,C), e_b(A,B,C), d_r(A,B,C),$ and~$d_b(A,B,C)$ in the natural way. When one of the sets $A,B$, or $C$ consists of a single vertex, we suppress the set brackets from the notation.

We sometimes write $A\sqcup B$ instead of $A\cup B$ when we want to emphasise that the sets $A$ and~$B$ are disjoint.

\paragraph{Organisation of the paper}
In \cref{sec:lowerbounds}, we present the constructions proving our lower bounds, namely \cref{prop:ell-paths,prop:loose_path_lower,prop:loose_cycle_lower,prop:non_transitive,prop:transitive_lb,prop:proj_planes}. The following three sections are devoted to our upper bound proofs. In \cref{sec:tools} we develop one of our main tools, which we call clique chains. Then we prove our upper bounds for loose paths and cycles in \cref{sec:loose_paths_up} and those for tight paths and cycles in \cref{sec:transitive_ub}. We end with a brief discussion of possible future directions in \cref{sec:conclusion}.

%%%%%%%%%%%%%%%%%%%%%%%%%%%%%%%%%%%%%%%%%%%
\section{Lower bound proofs}\label{sec:lowerbounds}

We begin by proving~\cref{{prop:ell-paths}}.

\begin{proof}[Proof of~\cref{prop:ell-paths}]
Let $\chi = \chi (H)$, $N = (\chi -1)(n-1) +\floor*{\frac{n}{k}}-1$, and $K = K^{(k)}_N$. Partition the vertex set of $K$ into $\chi$ sets $V_1,\dots, V_\chi$, where $\size{V_i} = n-1$ for all $i\in [\chi-1]$ and $\size{V_\chi} = \floor*{\frac{n}{k}}-1$. We then colour the edges of $K$ as follows. The red subgraph consists of all edges contained in a single $V_i$ for all $i\in [\chi]$ and those edges containing at least one vertex of $V_\chi$ and intersecting every other $V_i$ in at most $\ell-1$ vertices. All remaining edges are coloured blue. This colouring is illustrated in~\cref{fig:ell-paths}.
    
\begin{figure}[ht]
    \centering
    \includegraphics[scale=0.8,page=1]{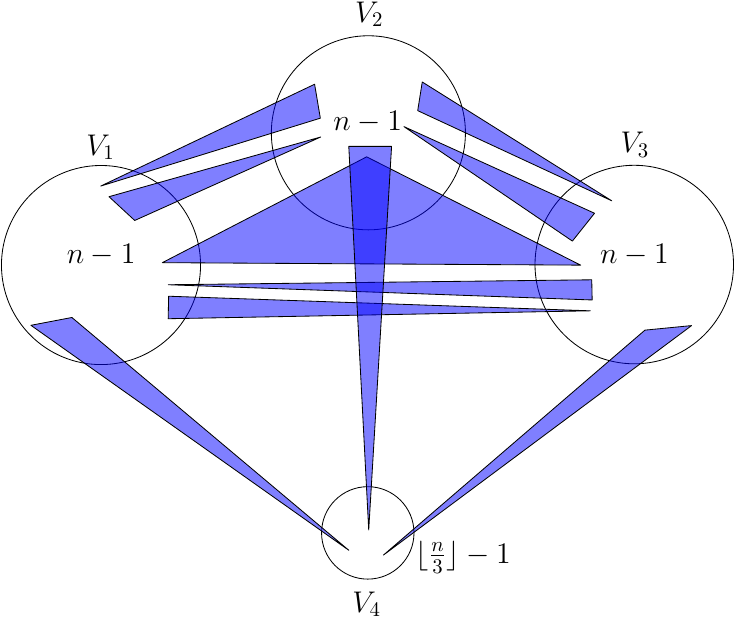}
    \caption{The blue subgraph in the proof of \cref{prop:ell-paths} for $k=3$, $\ell=2$, and $\chi=4$.}
    \label{fig:ell-paths}
\end{figure}

Suppose first that there exists a red copy $P$ of $P^{(k)}_{n,\ell}$. If $P$ contains an edge that is fully contained in $V_i$ for some $i\in [\chi-1]$, then  since $\ell\geq 2$, the entire path $P$ must be fully contained inside $V_i$. This implies that $v(P) \leq \size{V_i} = n-1$, a contradiction. Therefore, every edge of $P$ intersects at least two different sets $V_i$ for $i\in[\chi]$ or is fully contained in $V_\chi$. In particular, every edge of $P$ contains at least one vertex of $V_\chi$. Since $P$ contains a matching of size $\floor*{\frac{n}{k}}$, we have $\size{P\cap V_\chi}\geq \floor*{\frac{n}{k}} > \size{V_\chi}$, a contradiction.

Now suppose there exists a blue copy $H'$ of $H$ in $K$. This copy $H'$ has a naturally induced proper $\chi$-colouring with colour classes $W_i = V(H')\cap V_i$ for all $i\in [\chi]$. Now, in the induced subgraph $K[\bigcup_{i\in[\chi]}W_i]$ every blue edge $e$ intersecting $W_\chi$ intersects some set $W_j$ in at least $\ell$ vertices, a contradiction to our assumption about $H$.
\end{proof}

We now deduce \cref{prop:proj_planes} from \cref{prop:ell-paths}.
\begin{proof}[Proof of \cref{prop:proj_planes}]
    Consider any proper 2-colouring of $\cL_q$ with colour classes $V_1$ and $V_2$, which exists by Richardson~\cite{richardson1956finite}. Recall that $v(\cL_q) = q^2+q+1$, $k=q+1$, and every pair of vertices is contained in precisely one edge. Then, for each $i\in [2]$, the colour class $V_i$ has size at least two and there exists an edge $e$ such that $\size{e\cap V_i}\leq k-2$. The bound then follows from \cref{prop:ell-paths}. 
\end{proof}

\subsection{Loose paths and cycles}

We begin by presenting our construction for loose paths, after which we sketch the (very similar) argument for loose cycles. Recall that we define $\tau(k, \alpha)$ to be the largest $n$ such that there exists an $n$-vertex $k$-graph with independence number less than $\alpha$ and no two-edge loose path (see \cref{def:tau}).

\begin{proof}[Proof of \cref{prop:loose_path_lower}]
Let $N = (\chi-1)(n-1) + \tau(k-1,t) -2(k-1)$ and $K = K^{(k)}_N$. Partition the vertex set of $K$ into sets $V_1,\dots, V_\chi$, where $\size{V_i} = n-1$ for all $i\in [\chi-2]$, $\size{V_{\chi-1}} = n-2k+1$ and $\size{V_\chi} = \tau(k-1,t)$.
Let $J$ be a $(k-1)$-graph on $V_{\chi}$ with independence number less than~$t$ and no $(k-1)$-uniform loose path of length two, which exists by the definition of~$ \tau(k-1,t)$.
We colour an edge $e$ of $K$ red if $e$ is fully contained in  $V_i$ for some $i\in [\chi-1]$ or if $\size{e \cap V_{\chi-1}} = 1$ and $e \cap V_{\chi}$ is an edge in~$J$; the remaining edges are blue.
We now show that there is no red copy of~$P^{(k)}_{n,1}$ and no blue copy of~$H$ in $K$. The colouring is illustrated in \cref{fig:loose_path_lb}.

\begin{figure}[ht]
    \centering
    \includegraphics[scale=0.8,page=3]{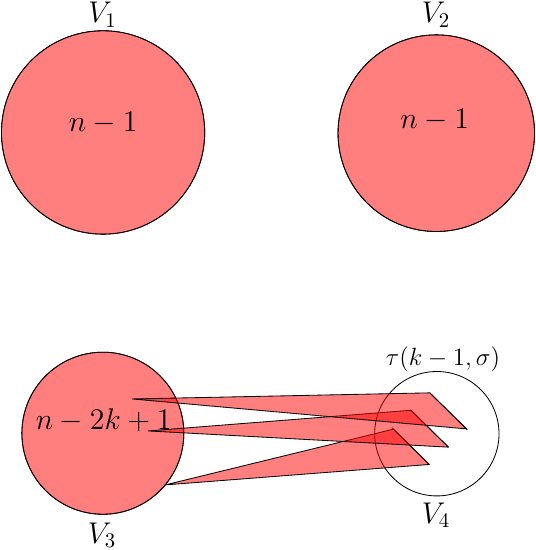}
    \caption{The red subgraph in the colouring from the proof of \cref{prop:loose_path_lower} for $\chi=4$.}
    \label{fig:loose_path_lb}
    \end{figure}

Consider a red loose path~$P$ in $K$ under this colouring. If $V(P) \cap V_i \ne \emptyset$ for some $i\in [\chi-2]$, then $V(P)$ is fully contained in $V_i$, so $v(P) < n$. So assume $V(P) \subseteq V_{\chi-1}\cup V_{\chi}$. If $e$ is an edge in~$P$ with $e \cap V_{\chi} \ne \emptyset$, then we claim that $e$ is either the first or the last edge of $P$. Indeed, note that $\size{e \cap V_{\chi}} = k-1$. If $e$ is not the first or last edge of $P$, there exists a red edge $e'$ in $P$  with $\emptyset \subsetneq e\cap e'\subseteq V_\chi$, and so
$\size{e' \cap V_{\chi}} = k-1$ and $\size{e\cap e'} = 1$. But then the $(k-1)$-edges $e \cap V_{\chi}$ and $e' \cap V_{\chi}$ form a loose path of length two in $J$, a contradiction. Thus $P' = P[V_{\chi-1}]$ is a loose path with $v(P')\geq v(P)-2(k-1)$. 
Hence $v(P) \leq v(P')+2(k-1) = n-2k+1+2k-2 = n-1$.

Suppose now that there is a blue copy of $H$ in $K$ with vertex classes $W_1,\dots, W_\chi$, where $W_\chi$ is the smallest class. By our assumption on the sizes of the $W_i$, for each $i\in [\chi-1]$, we have $\size{W_i\setminus V_\chi} > (\chi-1)(k-2)$, so there exists some $j_i\in [\chi-1]$ with $\size{W_i\cap V_{j_i}}\geq k-1$. Since the edges within a single $V_{j}$ are coloured red, we have $V_{j_i}\cap W_{i'} = \emptyset$ for all distinct $i,i'\in [\chi]$. Thus, all $j_i$ are distinct and we may assume that $j_i=i$ for all $i\in [\chi-1]$. It also follows that $W_\chi\subseteq V_\chi$ and thus $J[W_\chi]$ contains an edge. But all edges of $K$ intersecting $V_\chi$ in an edge of $J$ are red, a contradiction. 
\end{proof}

\begin{proof}[Proof of \cref{prop:loose_cycle_lower}]
The proof that $R(C_{n,1}^{(k)}, H)> (\chi-1)(n-1)+\tau(k-1,t)$ is very similar to that of \cref{prop:loose_path_lower}, except that now we have $\chi-1$ vertex classes of size $n-1$ and one of size $\tau(k-1,t)$. The rest of the argument is essentially the same.  

Assume now that $q$ is such that $\chi> \binom{q}{k-1}$. To show that $R(C_{n,1}^{(k)}, H)> (\chi-1)(n-1)+q$, we again let $N = (\chi-1)(n-1) + q$ and $K = K^{(k)}_N$. 
Partition the vertex set of $K$ into sets $V_1,\dots, V_\chi$, where $\size{V_i} = n-1$ for all $i\in [\chi-1]$ and $\size{V_\chi} = q$.
Write $\binom{V_{\chi}}{k-1} = \{S_i : 1\leq i \leq \binom{q}{k-1} \}$.
We colour an edge $e$ of $K$ red if $e$ is fully contained in some set $V_i$ for $i\in [\chi-1]$ or if $\size{e \cap V_{i}} = 1$ and $e \setminus  V_{i} = S_i$; the remaining edges are blue. 

The arguments showing that there is no monochromatic copy of $C_{n,1}^{(k)}$ in red and $H$ in blue are essentially the same as those in the proof of \cref{prop:loose_path_lower}.
\end{proof}

\subsection{Tight paths and cycles}

We now move on to our result concerning hypergraphs associated to non-transitive tournaments.
\begin{proof}[Proof of~\cref{prop:non_transitive}]
    Let~$N = (m-1)t$ and~$K = K^{(3)}_N$. We partition the vertex set of~$K$ into sets~$V_1,\dots, V_{m-1}$  satisfying~$\size{V_i} = t$ for all~$i\in [m-1]$. We then colour every $V_iV_iV_j$-edge for all~$1\leq i\leq j\leq m-1$ red and every other hyperedge blue. See~\cref{fig:non_transitive} for an illustration.

    \begin{figure}[ht]
    \centering
    \includegraphics[scale=0.8,page=2]{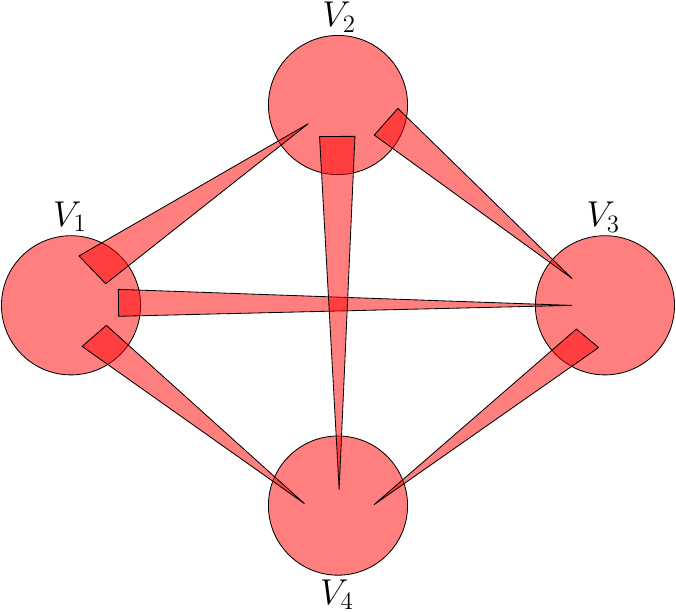}
    \caption{The red subgraph from the proof of \cref{prop:non_transitive} for $m=5$.}
    \label{fig:non_transitive}
    \end{figure}

    It is not difficult to see that a red tight path in this colouring has at most~$n-1$ vertices. Indeed, any red tight path must contain at most $b\leq t$ vertices from a single~$V_i$ plus at most $\floor*{\frac{b}{2}}+1$ vertices from~$V_{i+1} \cup\dots\cup V_{m-1}$, so its number of vertices cannot exceed~$t+\floor*{\frac{t}{2}}+1 < n$. 

    Now suppose there is a blue copy~$H'$ of~$H$ in~$K$ with vertex classes~$W_1,\dots, W_\chi$. For each {$j\in[\chi]$}, we have~$\size{W_j} = m$, and thus there exists an index~$k_j\in [m-1]$ such that~$\size{W_j\cap V_{k_j}}\geq 2$. Note that, since the edges fully contained in a single set~$V_i$ are red, for every arc~$(j,\ell)$ of~$T_{\chi}$, no set~$V_i$ can contain three vertices~$x,y,z$ such that~$x,y\in W_j$ and~$z\in W_\ell$. Therefore, all~$k_j$ are distinct.
    By the definition of colouring of $K$, it follows that $H'\brackets*{\bigcup\limits_{j\in [\chi]}\parens*{W_j\cap V_{k_j}}}$ is a tournament hypergraph  associated to the transitive tournament~$TT_{\chi}$, which contradicts the fact that~$H$ is a tournament hypergraph associated to a non-transitive tournament. 
\end{proof}

Finally, we prove our lower bound for tournament hypergraphs associated to transitive tournaments.
\begin{proof}[Proof of~\cref{prop:transitive_lb}]
     Let $R = \vec{R}(\chi)$ and~$T_{R-1}$ be a tournament on vertex set~$[R-1]$ that does not contain a copy of~$TT_{\chi}$, which exists by the definition of~$R$. Let~$N = \parens*{\floor*{\frac{2}{3}n}-2}(R-1)\geq$ $ \parens*{\frac{2}{3}n-3}(R-1)$ and~$K = K_N^{(3)}$. Partition the vertex set of~$K$ into sets~$V_1, \dots, V_{R-1}$ with~$\size{V_i} = \floor*{\frac{2}{3}n}-2$ for all~$i\in [R-1]$. We now assign the colour red to all $V_iV_iV_i$-edges for $i\in [R-1]$ and to all $V_iV_iV_j$-edges where~$(i,j)$ is an arc of~$T_{R-1}$. All remaining edges are coloured blue.

    Using a similar argument as in the proof of~\cref{prop:non_transitive}, we conclude that there is no red tight path on~$n$ vertices.  Suppose there exists a blue copy~$H'$ of~$H$ with vertex classes~$W_1,\dots, W_{\chi}$. Since for each~$j\in [\chi]$ we have~$\size{W_j}\geq R$, we know that there exists an integer~$k_j\in [R-1]$ such that~$\size{W_j\cap V_{k_j}}\geq 2$. As in the proof of \cref{prop:non_transitive}, all of these~$k_j$ are distinct.
    But then the hypergraph ~$H'\brackets*{\bigcup\limits_{j\in [\chi]} \parens*{W_j\cap V_{k_j}}}$ is a tournament hypergraph associated to~$TT_{\chi}$. But~$T_{R-1}$ does not contain a copy of~$TT_{\chi}$, a contradiction.
\end{proof}

\section{Tools}\label{sec:tools}

Most of this section is devoted to proving the existence of certain substructures that we call clique chains in $k$-graphs with a particular structure. Before we move on to clique chains, we state a useful fact and prove \cref{prop:tau_bounds}.

In order to apply the clique chain lemma that we will obtain later in this section, we will need a $k$-graph $G$ almost all of whose vertices can be partitioned into cliques of a certain size. We will find the required partition by using the following simple fact, shown by repeatedly applying the definition of Ramsey numbers.

\begin{fact}\label{lem:monochromatic_cliques}
Let~$k\geq 3$, $a,b\geq 1$, and $N\geq 1$ be integers. Then, for any red/blue-edge-colouring of~$K_N^{(k)}$, we can partition the vertex set of~$K_N^{(k)}$ into sets~$D_1,\dots, D_t, J$, where~$\size{J} < R(K_{a}^{(k)}, K_{b}^{(k)})$ and~$D_i$ induces a red copy of~$K_{a}^{(k)}$ or a blue copy of~$K_{b}^{(k)}$ for each~$i\in [t]$.  
\end{fact}

We also prove our bounds on $\tau(k,\alpha)$ in the non-trivial regime. The upper bound is due to Shagnik Das (personal communication). 

\begin{proof}[Proof of \cref{prop:tau_bounds}]
    We begin with the upper bound. Recall that a vertex cover of a hypergraph $G$ is a set $S\subseteq V(G)$ such that $e\cap S\neq \emptyset$ for every edge $e\in E(G)$. Let $G$ be a $k$-graph on $n$ vertices with independence number less than $\alpha$ and no loose path of length two. We will show that $n\leq 2\alpha-2$. Let $S$ be a minimal vertex cover of $G$. Then $V(G)\setminus S$ is an independent set and hence contains at most $\alpha-1$ vertices. By the minimality of $S$, for every vertex $v\in S$, there exists an edge $e$ of $G$ such that $e\cap S = \set{v}$. Therefore, $S$ cannot contain an edge either, so it also contains at most $\alpha-1$ vertices. Hence, $n\leq 2\alpha-2$, as required.

    For the lower bound, we write
    $\alpha - 1 = r(k-1)+s$, where $r = \floor*{\frac{\alpha-1}{k-1}}$ and $0\leq s\leq k-2$, and set $n = r(2k-2)+s = \alpha-1 + (k-1)\floor*{\frac{\alpha-1}{k-1}}$. Let $G$ be a $k$-graph on $n$ vertices consisting of $r$ copies of $K^{(k)}_{2k-2}$ and $s$ isolated vertices. Then $G$ contains no loose path of length two, and has independence number at most $r(k-1)+s = \alpha-1$.  Indeed, if a set $S$ contains more than $r(k-1)+s$ vertices, then it intersects one of the cliques in at least $k$ vertices, so $S$ cannot be independent.    
\end{proof}

\subsection{Clique chains}
Throughout this section, we assume that $k$ and $\ell$ are fixed and we might sometimes suppress~$\ell$ from the notation. In addition, we set $q_0 = q_0(k,\ell)$ to be the minimum length of an $\ell$-path containing at least $\max\set{k,2\ell}$ vertices. Note that $q_0  = \ceil*{\ell/(k-\ell)}$.

\begin{mydef}\label{def:clique_chain}
Let $k\geq 3$ and $1\leq \ell\leq k-1$ be integers. We call a $k$-graph~$F$ an \emph{open $k$-uniform {($\ell$-)}clique chain} if $V(F) = \{ v_i : i \in [p]\}$ and there exist intervals $I_1, \dots, I_d\subseteq [p]$ with $\bigcup_{j \in [d]} I_j = [p]$ such that
    \begin{enumerate}[label=(\alph*)]
        \item $\size{I_j} \ge k$ and $\size{I_j}\equiv \ell \pmod{k-\ell}$ for all $j\in [d]$;
        \item the set $\{v_i: i \in I_j\}$ induces a clique in $F$ for each $j\in [d]$;
        \item $\size{I_j \cap I_{j+1} } = \ell$ for all $j \in [d-1]$.
    \end{enumerate}
    Let $S_j = \{ v_i : i \in I_j\}$ for all $ j \in [d]$; we call the sets $S_j$ the \emph{elements} of $F$. 
    We say that $S_j$ is a \emph{flexible element  of~$F$} if $\size{S_j} >\max \{ k, 2 \ell\}$, and a \emph{rigid element  of~$F$} otherwise. 
    A vertex $v_i \in V(F)$ is a \emph{spine vertex} if $i \in I_{j} \cap I_{j+1}$ for some~$j$, and a \emph{flexible vertex} otherwise.

We define a \emph{closed $(\ell$-)clique chain} analogously by considering $V(F) = \{ v_i : i \in \mathbb{Z}/ p \mathbb{Z}\}$, i.e., the vertices of $F$  are ordered cyclically, and including the additional requirement that $\size{I_d \cap I_{1} } = \ell$. 
\end{mydef}
 
For simplicity, we will usually specify a clique chain by simply listing its elements, from which we can then deduce a suitable vertex ordering and choice of intervals. Note that if $\size{S_j} = k$ for all~$j \in [d]$, then $F$ is a $k$-uniform $\ell$-path (or $\ell$-cycle if $F$ is closed). More generally, we can think of an $\ell$-clique chain as a more flexible version of an $\ell$-path,
in particular because it contains multiple spanning $\ell$-paths. We will need this flexibility when we attempt to absorb extra vertices.
The following proposition is immediate from the definition. 

\begin{prop}\label{lem:chain_to_path}
Let $k\geq 3$ and $1\leq \ell\leq k-1$ be integers and $F$ be a $k$-uniform $\ell$-clique chain. 
Then $F$ has a spanning $\ell$-path if it is open, and a spanning $\ell$-cycle if it is closed. 
\end{prop}

In the rest of this section, we prove the following lemma, which guarantees the existence of $\ell$-clique chains satisfying additional properties in a certain class of hypergraphs. Before we state the lemma, we introduce a bit of additional terminology. 
Let $P$ be a $k$-uniform $\ell$-path. An \emph{end} of $P$ is the set of the first $\ell$ vertices or the last $\ell$ vertices of $P$.
For vertex sets $S$ and~$S'$, we say that $P$ \emph{connects $S$ and $S'$} if one end of~$P$ is contained in~$S$ and the other in~$S'$.

\begin{lem}\label{lem:alpha_chains}
    For all integers $k \ge 3$, $\ell \in \set*{1, k-1}$, $\alpha\geq 2$ and all $0<\eps<1$, there exists an integer $M_{\ref{lem:alpha_chains}}=M_{\ref{lem:alpha_chains}}(k,\alpha,\eps)$ such that, for all $M\geq M_{\ref{lem:alpha_chains}}$, the following holds.
		Let $G$ be a $k$-graph on vertex set $V_0 \sqcup V_1\sqcup \dots\sqcup V_t$, where $\size{V_i} = M$ and $V_i$ induces a copy of $K_{M}^{(k)}$ in $G$ for all $i\in [t]$.
    Suppose that, for any subset $I \subseteq [t]$ of size~$\alpha$ and any $W_i \subseteq V_i$ of size~$M/2$ for $i \in I$, the $k$-graph $G[ \bigcup_{i \in I} W_i ]$ contains an $\ell$-path of length~$q_0(k,\ell)$ connecting $W_j$ and~$W_{j'}$ for some distinct $j,j' \in I$.
		Then $G$ contains vertex-disjoint closed $\ell$-clique chains $Q_1, \dots, Q_{\alpha_0}$ with $\alpha_0 < \alpha$ satisfying: 
    \begin{enumerate}[label = \normalfont(\alph*)]
        \item $\size{  \bigcup_{i \in [t]}V_i \setminus \bigcup_{i \in [\alpha_0] } V( Q_{i} )} \le 10 t k$;\label{lem:alpha_chains:order}
        \item each chain $Q_i$ contains between $0.9\frac{ v(Q_i)}{M}$ and $1.1\frac{ v(Q_i) }{M}$ flexible elements; \label{lem:alpha_chains:num_flexible}
        \item for each $i\in [\alpha_0]$, $Q_i$ contains at least $(1- \eps)  v(Q_i)$ flexible vertices;\label{lem:alpha_chains:most_vxs}
        \item every flexible element $S$ satisfies $\size{S}\geq (1- \eps)M$;\label{lem:alpha_chains:flexible}
		\item for every pair of flexible elements $S,S'$ from different chains and every pair of vertex-disjoint $\ell$-paths $P_1,P_2$ connecting $S$ and $S'$ and satisfying $v(P_1),v(P_2)\leq 2k$, the set \mbox{$V(P_1)\cup V(P_2)$} contains a spine vertex of some chain. \label{lem:alpha_chains:in_between}
    \end{enumerate}

\end{lem}

Before we proceed with the proof of our main lemma, we introduce the notion of a path system. Intuitively, a path system is a collection of vertex-disjoint $\ell$-paths that serve as ``bridges'' between the different sets $V_i$ and allow us to connect several of these cliques into a single long $\ell$-path.

\begin{mydef}\label{def:path_system}
    Let $k\geq 3$ and  $1\leq \ell\leq k-1$ be integers. 
		Let $V = V_1\sqcup\dots\sqcup V_t$ be a set of vertices.
		Let~$G$ be a~$k$-graph on vertex set~$V$ and $F$ be a forest with $V(F) = [t]$.
    We say that a collection $\mathcal{L}$ is a \emph{$(k,\ell,F)$-path system} if: 
		\begin{enumerate}[label=(\alph*)]
        \item $\mathcal{L} = \{L_{e,1}, L_{e,2} : e \in E(F) \}$ consists of $2e(F)$ vertex-disjoint $k$-uniform $\ell$-paths;
        \item for every $e =xy \in E(F)$ and $i \in [2]$, the path $L_{e,i}$ connects $V_x$ and~$V_y$ and has order at most~$2k$.			\end{enumerate}
        We sometimes consider $\mathcal{L}$ as a $k$-graph with vertex set $V(\mathcal{L}) = \bigcup_{L\in \mathcal{L}} V(L)$ and $ E(\mathcal{L}) = \bigcup_{L\in \mathcal{L}} E(L)$.
\end{mydef}

The first step towards proving \cref{lem:alpha_chains} is to establish that, under similar assumptions, we can find a path system with fewer than $\alpha$ components that does not use too many vertices from any one clique. Once we have shown \cref{lem:path_system} below, all that is left to do to prove \cref{lem:alpha_chains} is to show that going from clique to clique along the bridges given by the path system will produce the required clique chains.

\begin{lem}\label{lem:path_system}
    For all integers $k \ge 3$, $\ell \in \set*{1, k-1}$, $\alpha\geq 2$ and all $0<\eps<1$, there exists an integer $M_{\ref{lem:path_system}} = M_{\ref{lem:path_system}}(k,\alpha,\eps)$ such that, for all $M\geq M_{\ref{lem:path_system}}$, the following holds. 
		Let $G$ be a $k$-graph on vertex set $V_0 \sqcup V_1\sqcup \dots\sqcup V_t$, where $\size{V_i} = M$ and $V_i$ induces a copy of $K_{M}^{(k)}$ in $G$ for all $i\in [t]$.
    Suppose that, for any subset $I \subseteq [t]$ of size~$\alpha$ and any $W_i \subseteq V_i$ of size~$M/2$ for $i \in I$, the $k$-graph $G[ \bigcup_{i \in I} W_i ]$ contains an $\ell$-path of length $q_0(k, \ell)$ connecting $W_j$ and~$W_{j'}$ for some distinct $j,j' \in I$.
	Then $G$ contains a $(k,\ell,F)$-path system~$\mathcal{L}$ such that 
	\begin{enumerate}[label={\rm(\alph*)}]
		\item $F$ is some forest with $V(F) = [t]$ containing fewer than $\alpha$ components;
		\item $\size{ V(\mathcal{L}) \cap V_i } \le \eps M$ for all $i \in [t]$;
		\item for every component $T$ of~$F$, we have $ \size{ \bigcup_{e \in E(T)} V(L_{e,1} \cup L_{e,2}) 	\setminus  \bigcup_{i \in V(T)} V_i } \le 4 \alpha k$; 
		\item for any $i$ and $i'$ contained in different components of $F$, the $k$-graph $G \setminus V( \mathcal{L}) $ does not contain two vertex-disjoint $\ell$-paths, each of order at most~$2k$, connecting~$V_i$ and~$V_{i'}$.\label{lem:path_system:no_short_paths}
	\end{enumerate}
\end{lem}

\begin{proof}
Our proof strategy is as follows: starting with a  $(k,\ell, F)$-path system $\cL$ such that $F$ contains too many components, we will show that each component $T_j$ contains a vertex $i_j$ such that $\cL$ uses very few vertices from $V_{i_j}$; as a result, we will be able to add new edges to $F$ connecting some of these $i_j$ while maintaining the property that the path system intersects each set $V_i$ in few vertices. Finally, we will add some extra edges to $F$ to ensure that property~\ref{lem:path_system:no_short_paths} holds. We now present the formal proof.

Set 
\begin{align*}
M_{\ref{lem:path_system}} = \max\set{20k\alpha/\eps, 10\alpha^2k, 12k^2 \alpha^k}
\end{align*}
and let $M\geq M_{\ref{lem:path_system}}$.
Suppose that $G$ satisfies the hypothesis. 
We first prove that $G$ contains a $(k,\ell,F)$-path system~$\mathcal{L}$ for some forest $F$ with $V(F) = [t]$ containing fewer than $\alpha$ components such that
\begin{enumerate}[label=(\roman*)]
    \item $\size{ V(\mathcal{L}) \cap V_i } \le 2\eps M/3$ for all $i \in [t]$;\label{item:intersection}
    \item for every component $T$ of~$F$, we have $\bigcup_{e \in E(T)} V(L_{e,1} \cup L_{e,2}) 	\subseteq \bigcup_{i \in V(T)}V_i$. \label{item:containment}	
\end{enumerate}
We will then modify this path system to prove the properties required in the statement.

First note that the empty path system is a $(k,\ell,\overline{K_t})$-path system satisfying~\ref{item:intersection} and~\ref{item:containment}, where~$\overline{K_t}$ denotes the empty graph on $[t]$. Now, let $\cL$ be a  $(k,\ell,F)$-path system satisfying~\ref{item:intersection} and~\ref{item:containment}, where $V(F) = [t]$ and the number of components in $F$ is minimal.

Suppose that $F$ contains $\alpha$ components, call them $T_1, \dots, T_{\alpha}$.
Consider $j \in [\alpha]$. 
Note that $\size{ V(\mathcal{L}) \cap \bigcup_{i \in V(T_j)}V_i } \le 4k v(T_j)$, since $\mathcal{L}$ contains at most $2v(T_j)$ paths, each of order at most $2k$,  associated with $T_j$.
If $\size{ V(\mathcal{L}) \cap V_i} \geq  \eps M/3$ for all $i\in V(T_j)$, we then have 
\begin{align*}
4 k v(T_j)\geq \size*{ V(\mathcal{L}) \cap \bigcup_{i \in V(T_j)}V_i } \geq \frac{\eps M v(T_j)}{3}.
\end{align*}
This implies that $M \le 12 k /\eps$, a contradiction.
Thus, there exists an $i_j \in V(T_j)$ such that $\size{ V(\mathcal{L}) \cap V_{i_j} } \le \eps M/3$.

\medskip
Let $I = \{i_j : j \in [\alpha]\}$ and $W_i =  V_{i} \setminus V(\mathcal{L})$ for all $i\in I$, and notice that $\size{W_i} \geq \parens*{1- \eps/3}M$.

\begin{clm}\label{clm:pathsystem}
The $k$-graph $G[ \bigcup_{i \in I} W_i ]$ contains a $(k,\ell,F')$-path system $\mathcal{L}'$ for some nonempty forest~$F'$.
Moreover, $V(\mathcal{L}') \subseteq \bigcup_{i \in V(F')}V_i$ and $v(\mathcal{L}') \le 4k^2$.
\end{clm}

\begin{proof}[Proof of \cref{clm:pathsystem}]
    We first prove the case when $\ell = k-1$. 
    Consider a maximal collection~$\cL'$ of vertex-disjoint $(k-1)$-paths of length~$q_0$ in~$G[ \bigcup_{i \in I} W_i ]$, each of which connects $W_j$ and $W_j'$ for distinct $j,j'\in I$. 
    We claim that $\size{\cL'} > \binom{\alpha}{2}$. 
    Suppose otherwise, and note that $v(\cL')  \leq \binom{\alpha}{2}2k$.
		For each $i\in I$, let $W'_i = W_i\setminus V(\cL')$, so $\size{W'_i}\geq (1-\eps/3) M - \alpha^2 k \geq M/2$.
		Recall that $|I| = \alpha$, so by assumption there exists a $(k-1)$-path of length~$q_0$ connecting $W'_j\subseteq W_j$ and $W'_{j'}\subseteq W_{j'}$ for distinct $j,j'\in I$, contradicting the maximality of $\cL'$.
    
    Hence $\size{\cL'} > \binom{\alpha}{2}$ and thus there are two $(k-1)$-paths in~$\cL'$ connecting the same pair $W_j$ and~$W_{j'}$.
    Since each such path~$P$ has order~$2(k-1)$, we deduce that $V(P) \subseteq W_j \cup W_{j'}$.
    Then letting $F'$ be a forest on $[t]$ with $E(F') = \set{jj'}$ gives the desired path system. 
    \smallskip

    Now consider the case where $\ell = 1$. 
    By a similar argument as above, we can find a matching~$\mathcal{M}$ in~$G[ \bigcup_{i \in I} W_i ]$ of size at least $2k \alpha^k$ such that no edge lies entirely in some~$V_i$. 
    Thus we can find a submatching $\mathcal{M}' \subseteq \mathcal{M}$ of size at least~$2k$ such that $\size{ e \cap V_i } = \size{ e' \cap V_i }$ for all $i \in I$ and all $ e,e' \in \mathcal{M}'$. 
    Let $I' = \{ i \in I : V_i \cap V(\mathcal{M}') \ne \emptyset \}$. 
    Note that $2 \le \size{I'} \le k$ and any edge of $\mathcal{M}'$ is a $1$-path of length~$q_0 = 1$ connecting $V_j$ and $V_j'$ for all $j,j' \in I'$. 
    By removing edges from $\mathcal{M}'$ if necessary, we have $\size{ \mathcal{M}' }= 2 \size{I'}$. 
    Set $F'$ to be a forest on $[t]$ with precisely one nonempty tree on vertex set~$I$.
    Hence $\mathcal{M'}$ is a $(k,\ell,F')$-path system as desired.     
\end{proof}

Let $\mathcal{L}'$ and $F'$ be as given by~\cref{clm:pathsystem}.
Then $\mathcal{L} \cup \mathcal{L}'$ is a $(k,\ell,F \cup F')$-path system. 
Note that 
\begin{align*}
 \size{ V(\mathcal{L} \cup \mathcal{L'}) \cap V_i } \le 
\begin{cases}
	\eps M/3 + 4k^2 \leq 2 \eps M/3 & \text{ for all $i \in I$,}\\
	\size{ V(\mathcal{L}) \cap V_i } \le 2 \eps M/3 & \text{ for all $i \in [t] \setminus I$.}
\end{cases}
\end{align*}
This contradicts the choice of $\cL$, more precisely the fact that~$F$ minimises the number of components. 
Hence we may assume that $\mathcal{L}$ satisfies \ref{item:intersection} and~\ref{item:containment} and $F$ has fewer than $\alpha$ components. 

\medskip
We now complete the proof of the lemma by modifying~$\mathcal{L}$. 
Suppose that there exist $i$ and $i'$ in different components of~$F$ such that $G \setminus V( \mathcal{L}) $ contains two vertex-disjoint $\ell$-paths~$P$ and~$P'$, each of order at most~$2k$, connecting~$V_i$ and~$V_{i'}$. 
We then add the edge $ii'$ to~$F$ and $\set{P,P'}$ to~$\mathcal{L}$. 
Repeating this procedure at most $\alpha$ times gives a path system $\cL$ with $\size{V(\cL)\cap V_i} \leq 2\eps M/3 + 4k\alpha \leq \eps M$ for all $i\in [t]$ satisfying the required properties. 
Furthermore, for every component $T$ of~$F$, we have $ \size{ \bigcup_{e \in E(T)} V(L_{e,1} \cup L_{e,2}) 	\setminus  \bigcup_{i \in V(T)} V_i } \le 4 \alpha k$. 
\end{proof}

We are now ready to complete the proof of \cref{lem:alpha_chains}.

\begin{proof}[Proof of Lemma~$\ref{lem:alpha_chains}$]
As explained earlier, our goal now is to show that we can use the $(k,\ell, F)$-path system guaranteed by \cref{lem:path_system} to join the cliques into clique chains. For each component~$T$ of $F$, we will first build a closed walk in $T$ visiting every vertex of $T$, which will serve as a ``template'' for our clique chain. The edges of this walk will be replaced by paths from the path system. On the other hand, each vertex of the walk will become either a flexible element or a short $\ell$-path connecting two paths from the path system; in either case, the structure will use only vertices from the corresponding set $V_i$. We proceed with the details.

Without loss of generality, assume $\eps<1/10$. 
Let 
\begin{align*}
M_{\ref{lem:alpha_chains}} = \max\set{ M_{\ref{lem:path_system}}( k,\alpha, {\eps}/{(20k)}), 20\alpha k/\eps} 
\end{align*}
 and $M\geq M_{\ref{lem:alpha_chains}}$. 
Suppose that $G$ satisfies the hypothesis. 
We begin by applying \cref{lem:path_system} (with $\eps/(20k)$ playing the role of $\eps$) to obtain a $(k,\ell,F)$-path system~$\cL$ such that 
\begin{enumerate}[label={\rm(\alph*$'$)}]
		\item $F$ is some forest with $V(F) = [t]$ containing fewer than $\alpha$ components;
		\item $\size{ V(\mathcal{L}) \cap V_i } \le \eps M/(20k)$ for all $i \in [t]$;
		\item for every component $T$ of~$F$, we have $ \size{ \bigcup_{e \in E(T)} V(L_{e,1} \cup L_{e,2}) 	\setminus  \bigcup_{i \in V(T)} V_i } \le 4 \alpha k$;
		\item for any $i$ and $i'$ contained in different components of $F$, the $k$-graph $G \setminus V( \mathcal{L}) $ does not contain two vertex-disjoint $\ell$-paths, each of order at most~$2k$, connecting~$V_i$ and~$V_{i'}$. \label{itm:alpha:4}
\end{enumerate}

Now let $T$ be a component of $F$. We show how to build a closed $\ell$-clique chain from~$T$.
Let $w_1\dots w_{b} w_{1}$ be a closed walk in~$T$ that visits every vertex at least once and traverses every edge exactly twice.\footnote{To see that such a walk exists, consider~$T$ and build an auxiliary multigraph~$T'$ on the same vertex set in which, for every edge~$xy\in T$, we add two edges between~$x$ and~$y$ to~$T'$. Then~$T'$ has only even degrees and thus it has an Euler circuit.} 
Let $L_1,\dots, L_{b}$ be an ordering of the elements of $\set{L_{e,1}, L_{e,2}:e\in E(T)}$, where $L_j \in \set{L_{w_jw_{j+1},1}, L_{w_jw_{j+1},2}}$ for each $j\in [b]$, where we take $w_{b+1} = w_1$. 
We treat $L_j$ as an $\ell$-path from $V_{w_j}$ to~$V_{w_{j+1}}$.

We construct the required clique chain in two steps.  Recall that an $\ell$-cycle is itself an $\ell$-clique chain in which each clique consists of a single edge. 
In the first step, we join $L_1, \dots, L_{b}$ into an $\ell$-cycle as follows. Consider $j \in [b]$.
Let $x_1, \dots, x_{\ell}$ be the last $\ell$ vertices of~$L_{j-1}$ and $y_1, \dots, y_{\ell}$ be the first $\ell$ vertices of~$L_{j}$ (for convenience we set $L_0 = L_{b}$).
Thus, $x_1 \dots x_{\ell}, y_1, \dots, y_{\ell} \in V_{w_j}$.
We now reserve an arbitrary $\ell$-path~$L'_j$ of length $q_0$ starting with $x_1 \dots x_{\ell}$ and ending with $y_1, \dots, y_{\ell}$ such that $V(L_j') \subseteq V_{w_j}$.
Recall that each $V_i$ induces a copy of $K^{(k)}_M$ in~$G$, so $L_j'$ exists. 
For each $i \in [t]$, at most $\size{ V(\mathcal{L}) \cap V_i } \le \eps M/(20k)$ of the vertices~$w_j$ are equal to~$i$. 
Hence $\cL \cup \bigcup_{j \in [b]} L'_j$ together use up at most ${\eps}M/6$ vertices from each~$V_i$. 
Furthermore, we can assume that $L'_1, \dots,L'_{b}$ are vertex-disjoint. 
Therefore the concatenation of $L'_1,L_1,L'_2,\dots, L'_{b},L_{b}$ yields an $\ell$-cycle, which is itself a closed $\ell$-clique chain.

It remains to create the required number of flexible elements.  
For each $i \in V(T)$, pick an index $j\in [b]$ such that $w_j = i$, and set $U_i = V_i \setminus \parens*{V(\mathcal{L}) \cup \bigcup_{j \in [b]}V( L'_j)}$.
Recall that $V(L'_j) \subseteq V_{w_j} = V_i$.
Replace $V(L'_j)$ with a set~$W_i$, where $W_i$ is a subset of largest size satisfying $\size{W_i}\equiv \ell \pmod{k-\ell}$ and $ V(L'_j) \subseteq W_i \subseteq V(L'_j)\cup U_i$.
Note that 
\begin{align*}
\size{W_i} \geq (1- \eps/6) M-k\geq (1-\eps/3) M.
\end{align*}
This creates a closed $\ell$-clique chain, whose flexible elements are the~$W_i$ for $i \in V(T)$.
Denote the resulting clique chain by~$Q$. Then $Q$ satisfies property~\ref{lem:alpha_chains:flexible}. 
The number of flexible elements in~$Q$ is $v(T)$. 
Note that 
\begin{align*}
      \left|  \bigcup_{i \in V(T)} W_i \right| \le 
v(Q)  & \le  \left|  \bigcup_{i \in V(T)} V_i \right| + \left| \bigcup_{e \in E(T)} V(L_{e,1} \cup L_{e,2}) 	\setminus  \bigcup_{i \in V(T)} V_i \right|,
\end{align*}
or in other words
\begin{align*}
    (1-\eps/3) M v(T)  \le v(Q)&\le M v(T) + 4 \alpha k \le (1+\eps/3) M v(T).
\end{align*}
A routine calculation then shows that the number of flexible elements is between $0.9\frac{v(Q)}{M}$ and $1.1\frac{v(Q)}{M}$, showing that $Q$ satisfies property~\ref{lem:alpha_chains:num_flexible}.
Also, the number of flexible vertices in~$Q$ is at least 
\begin{align*}
    \sum_{i \in V(T)} ( \left|   W_i \right|-2k) \ge (1- 2 \eps/3) M v(T) \ge (1-\eps)v(Q).
\end{align*} 
Finally, note that $\size{ V_i \setminus V(Q \cup \mathcal{L}) } \le k$ for all $i \in V(T)$.

Repeating this for every component of~$F$, we obtain  (closed) $\ell$-clique chains $Q_1,\dots, Q_{\alpha_0}$, where $\alpha_0$ is the number of components of~$F$. 
Note that \ref{lem:alpha_chains:order} holds by the final observation above.  
Parts~\ref{lem:alpha_chains:num_flexible},~\ref{lem:alpha_chains:most_vxs}, and~\ref{lem:alpha_chains:flexible} are verified above.
Finally, \ref{lem:alpha_chains:in_between} follows from~\ref{itm:alpha:4}.
\end{proof}

%%%%%%%%%%%%%%%%%%%%%%%%%%%%%%%%%%%%%%%%%%%

\section{Upper bounds for loose paths and cycles}\label{sec:loose_paths_up}

Throughout this section, $H$ will be a fixed $k$-graph for some $k\geq 3$, and we will write $m = v(H)$, $\chi = \chi(H)$, and $\sigma = \sigma(H)$. 
As the proofs of \cref{thm:loosepaths_upper,thm:loosecycles_upper} are very similar,
we only present the proof of \cref{thm:loosepaths_upper} and give a sketch 
 to the proof of \cref{thm:loosecycles_upper}. We will make use of the following simple fact.

\begin{fact}\label{fact:independence_nr_F}
Let $m, \chi, k\geq 1$ be integers with $k\geq 3$.
Suppose that the complete $k$-graph $K^{(k)}_{\chi m}$ is red/blue-edge-coloured.
Let $V(K^{(k)}_{\chi m})$ be partitioned into subsets $W_1,\dots, W_{\chi}$, each of size~$m$.
Then at least one of the following holds:
\begin{enumerate}[label={\rm(\alph*)}]
    \item There exists a red edge in $K^{(k)}_{\chi m}$ connecting $W_j$ and $W_{j'}$ for some distinct $j,j'\in [\chi]$.\label{lem:independence_nr_F:red}
    \item Every edge in $\bigcup_{i \in [\chi]}W_i$ that is not entirely contained in some $W_i$ is blue. In particular, the $k$-graph induced by $\bigcup_{i \in [\chi]}W_i$ contains a blue copy of $H$.\label{lem:independence_nr_F:blue}
\end{enumerate}
\end{fact}

We are now ready to derive a statement guaranteeing the existence of a suitable collection of red 1-clique chains.

\begin{cor} \label{clm:chainspartition}
Let $\chi, k \geq 1$ be integers with $k\geq 3$.
Then there exists an integer $M_{\ref{clm:chainspartition}} = M_{\ref{clm:chainspartition}}(k,\chi)$ such that, for all $M \ge M_{\ref{clm:chainspartition}}$ and $m\geq 1$, there is an integer
$n_{\ref{clm:chainspartition}} = n_{\ref{clm:chainspartition}}(k,\chi,m, M)$ for which the following holds. Let $n \ge n_{\ref{clm:chainspartition}}$ and $N \ge (\chi-1)(n-1)$, and let $H$ be a $k$-graph with $v(H) = m$, $\chi(H) = \chi$, and $\sigma (H) = \sigma$.
Suppose $K = K^{(k)}_N$ is given a $(P_{n,1}^{(k)},H)$-free edge-colouring.
Then there exist vertex-disjoint red open $1$-clique chains $Q_1,\dots, Q_{\chi-1}$  such that
        \begin{enumerate}[label={\rm (\alph*)}]
					\item 
					\label{itm:clm:chainspartition:vertices}
					for all $i \in [\chi - 1] $, $v(Q_i) = n - \ell_i$ with $ \ell_i  \le  \frac{n}{100k}$ and $\ell_i \equiv 0 \pmod{k-1}$;
      		\item each chain $Q_i$ contains at most $1.2 n/ M $ flexible elements;\label{itm:clm:chainspartition:flexiblevxssize}
					    \item each chain $Q_i$ contains at least  $2n/3 + 10k \ell_i$ flexible vertices; \label{itm:clm:chainspartition:flexiblevxs}
					
					\item every flexible element contains at least $M/4$ vertices; 
					\label{itm:clm:chainspartition:flexibleelts}
        \item for each $i\in [\chi-1]$, the first and last element of $Q_i$ is flexible;\label{itm:clm:chainspartition:firstlast}
          \item 
					\label{itm:clm:chainspartition:crossing} every red loose path of length at most two connecting flexible elements from two different chains in $K$ contains a spine vertex;
					\item 
					\label{itm:clm:chainspartition:red}
					for every subset $B\subseteq V(K)$ of size $\sigma$  and every collection of subsets of flexible vertices $A_1,\dots, A_{\chi-1}$ with $A_i\subseteq V(Q_i)\setminus B$ and $\size{A_i}\geq m$ for all $i\in[\chi-1]$, there exists a red edge $e\subseteq A_\ell\cup B$ for some $\ell\in [\chi-1]$ satisfying $1\leq \size{e\cap B}\leq k-1$.
				\end{enumerate}
\end{cor}

\begin{proof}
 Let $\eps < \frac{1}{100k}$ be arbitrarily small and $M_{\ref{clm:chainspartition}} = \max\set*{\ceil*{50 k \chi/\eps}, M_{\ref{lem:alpha_chains}}(k,\alpha,\eps)}$.
Let ${M \ge M_{\ref{clm:chainspartition}}}$ and $m\geq 1$. Set
{$
n_{\ref{clm:chainspartition}} = \lceil 2 R(K_{M}^{(k)}, K_{m}^{(k)}) / \eps \rceil$} and let $n\geq n_{\ref{clm:chainspartition}}$. 
Suppose that $H,\sigma,N,K$ satisfy the hypothesis. 
Note that $K$ does not contain a blue copy of~$H$. 
Apply \cref{lem:monochromatic_cliques} to partition $V(K)$ into subsets $V_0, V_1, \dots,V_t$, 
where 
\begin{align}
	\size{V_0} < R(K_{M}^{(k)}, K_{m}^{(k)}) \leq \eps n/2, \label{eqn:|V_0|1}
\end{align}
each $V_i$ induces a red copy of $K_{M}^{(k)}$, and 
\begin{align*}
\frac{N-\eps n/2}{M}\leq t\leq \frac{N}{M} \leq \frac{\eps N}{50k\chi}. 
\end{align*}

Now, \cref{fact:independence_nr_F} allows us to apply \cref{lem:alpha_chains} (with $(\ell, \alpha) = (1, \chi)$ and the subgraph of $K$ induced by the red edges playing the role of $G$) and obtain vertex-disjoint red closed $1$-clique chains $Q_1, \dots, Q_{\alpha_0}$ with $\alpha_0 < \chi$ satisfying: 
    \begin{enumerate}[label = \normalfont(\alph*$'$)]
        \item $\size{  \bigcup_{i \in [t]}V_i \setminus \bigcup_{i \in [\alpha_0] } V( Q_{i} )} \le 10 t k$;\label{itm:cor:1}
        \item each chain $Q_i$ contains between $0.9\frac{ v(Q_i)}{M}$ and $1.1\frac{ v(Q_i) }{M}$ flexible elements; \label{itm:cor:2}
        \item for each $i\in [\alpha_0]$, $Q_i$ contains at least $(1- \eps)  v(Q_i)$ flexible vertices;\label{itm:cor:3}
        \item every flexible element $S$ satisfies $\size{S}\geq (1-\eps)M$;\label{itm:cor:4}
		\item for every pair of flexible elements $S,S'$ from different chains and every pair of vertex-disjoint $\ell$-paths $P_1,P_2$ connecting $S$ and $S'$ and satisfying $v(P_1),v(P_2)\leq 2k$, the set $V(P_1)\cup V(P_2)$ contains a spine vertex of some chain. \label{itm:cor:5}
		
    \end{enumerate}

Note that each $Q_i$ contains a red spanning $1$-cycle by \cref{lem:chain_to_path}, so $Q_i$ contains a red loose path on $v(Q_i)-(k-2)$ vertices. Hence $v(Q_i) < n+k-2$, and so  $v(Q_i) < n $. 
Together with  \eqref{eqn:|V_0|1} and \ref{itm:cor:1}, this implies
\begin{align*}
	N & \le \size{V_0} + \sum_{i \in [\alpha_0]} v(Q_i) + \left|   \bigcup_{i \in [t]}V_i \setminus \bigcup_{i \in [\alpha_0] } V( Q_{i} )  \right| 
	 \le \frac{\eps n}{2} + (\chi - 1)n + 10 t k \leq \frac{\eps n}{2} + (\chi - 1)n + \frac{\eps N}{5\chi}. 
\end{align*}
Rearranging, and solving for $N$, we obtain 
\begin{align}
    N\leq \frac{\chi-1+\eps/2}{1-\eps/(5\chi)}n \leq (\chi-1+\eps)n,\label{eq:boundN}
\end{align}
where the last inequality holds for all $0<\eps<1$.
\cref{prop:loosepaths_asymptotic} follows directly from \eqref{eq:boundN}. In addition, \eqref{eq:boundN} implies that
\begin{align}
\size*{ V(K) \setminus \bigcup_{i \in [\alpha_0] } V( Q_{i} ) } \le \frac{\eps n}{2} + 10tk \leq \frac{\eps n}{2} + \frac{\eps N}{5\chi} \leq \eps n .  \label{eqn:cor1}
\end{align} 
Hence we have $\alpha_0 = \chi - 1$ and 
\begin{align}
	v(Q_i) \ge (1- \eps) n .\label{eqn:cor2}
\end{align}

Suppose that there exist distinct $i,j \in [\chi-1]$ and a  red loose path~$P$ of length at most two connecting flexible elements from $Q_i$ and $Q_j$ that does not contain any spine vertices. 
Since~$Q_i$ and $Q_j$ are closed, there exists a red open clique chain~$Q$ in $K[V(Q_i) \cup V(Q_j) \cup V(P)]$ such that $\size{ (V(Q_i) \cup V(Q_j)) \setminus V(Q)} \le 4k$. 
By \cref{lem:chain_to_path}, there is a red loose path on $v(Q) \ge v(Q_i \cup Q_j)-4k > n$ vertices, a contradiction. 
Hence $Q_1, \dots, Q_{\chi-1}$ satisfy~\ref{itm:clm:chainspartition:crossing}.

We now `cut each $Q_i$ open' by simply splitting one of its flexible elements into two sets, and possibly discarding at most $2k$ vertices to adjust the parity, to obtain two flexible elements with at least $M/4$ vertices each. All other flexible elements remain intact and thus have size at least~$(1-\eps)M$.

We now verify that these $Q_i$ satisfy the desired properties.
Clearly properties \ref{itm:clm:chainspartition:flexibleelts} and \ref{itm:clm:chainspartition:firstlast}  hold by our construction and \ref{itm:clm:chainspartition:crossing} remains true. 
Note that~\eqref{eqn:cor2} implies~\ref{itm:clm:chainspartition:vertices}, and \ref{itm:clm:chainspartition:flexiblevxssize} follows from~\ref{itm:cor:2}. To see that part~\ref{itm:clm:chainspartition:flexiblevxs} holds,
note that properties~\ref{itm:clm:chainspartition:vertices} and~\ref{itm:cor:3} imply that the number of flexible vertices in each $Q_i$ is at least $(1-\eps)(n-\ell_i)$, which exceeds $2/3n+10k\ell_i$ if $\ell_i\leq n/(100k)$ and $\eps < 1/10$.

Finally, we now show~\ref{itm:clm:chainspartition:red}.
Let $B, A_1,\dots, A_{\chi-1}$ be as given as in~\ref{itm:clm:chainspartition:red}.
By~\ref{itm:clm:chainspartition:crossing}, all edges intersecting at least two different sets $A_i$ are blue.
If additionally all edges  $e\subseteq A_i\cup B$ satisfying $1\leq \size{e\cap B}\leq k-1$ are blue, then there exists a blue copy of $H$ in $K[B\cup \bigcup_{i \in[\chi-1]} A_i]$, a contradiction.
\end{proof}

Note that the proof of \cref{clm:chainspartition} (up to \eqref{eq:boundN}) readily yields \cref{prop:loosepaths_asymptotic}. 
We now proceed to the proof of \cref{thm:loosepaths_upper}.

\begin{proof}[Proof of \cref{thm:loosepaths_upper}]  
Before we begin, we briefly outline our strategy. We suppose we have a red/blue-colouring of a large complete $k$-graph $K_N^{(k)}$ with no red loose path on $n$ vertices and no blue copy of $H$. We start by finding $\chi-1$ red clique chains $Q_1,\dots, Q_{\chi-1}$, satisfying the properties given in \cref{clm:chainspartition} and maximising the total number of vertices among all such collections of chains, and a leftover set of vertices $W$. As a first step, we show that each chain~$Q_i$ contains at least one large flexible element $S_i$. As there is no blue copy of $H$, there are many red edges within the sets $W\cup S_i$. If there is a large red matching in some $W\cup S_i$ such that each edge in the matching intersects $S_i$ in at least two vertices, then we can extend the clique chain~$Q_i$ while preserving its properties, contradicting the maximality of our collection. We can then deduce that there are many red edges in each $W\cup S_i$ intersecting $S_i$ in precisely one vertex. Then we have two cases depending on how large $W$ is. If $W$ contains many vertices, then by considering the $(k-1)$-uniform graph on $W$ defined by these red edges, we find a two-edge $(k-1)$-uniform loose path that again allows us to extend some chain~$Q_i$. On the other hand, if $W$ is small, then we argue that each of the chains contains almost $n$ vertices and we can create a loose path on $n$ vertices by adding just one or two edges to some chain; doing so will yield our final contradiction. 

Set 
\begin{align*}
M =  \max\set*{M_{\ref{clm:chainspartition}}(k,\chi), 100\chi k^2, 20m +  4\binom{\tau(k-1,\sigma)+1}{k-1}} 
\text{ and } 
	n_0 =  n_{\ref{clm:chainspartition}}(k,\chi,m,M).
\end{align*}
Let 
\begin{align*}
n \ge n_0 \text{ with $n \equiv 1\pmod{k-1}$,} \quad
c = \max \{ \tau(k-1,\sigma) -2k+3, \sigma \}, \quad
N = (\chi-1)(n-1) +  c,
\end{align*}
and $K = K_N^{(k)}$.
Suppose to the contrary that $K$ has a $(P_{n,1}^{(k)}, H)$-free edge-colouring, and we fix one such edge-colouring.
Apply \cref{clm:chainspartition} and obtain vertex-disjoint red open 1-clique chains  $Q_1, \dots, Q_{\chi-1}$ such that 
        \begin{enumerate}[label={\rm (\alph*)}]
					\item 
					\label{itm:loosepath1}
					for all $i \in [\chi - 1] $, $v(Q_i) = n - \ell_i$ with $ \ell_i  \le  \frac{n}{100k}$ and $\ell_i \equiv 0 \pmod{k-1}$;
      		\item each chain $Q_i$ contains at most $1.2 n/ M $ flexible elements;
					\label{itm:loosepath2}
					\item each chain $Q_i$ contains at least  $2n/3 + 10k \ell_i$ flexible vertices; 
					\label{itm:loosepath3}
					\item every flexible element contains at least $M/4$ vertices; 
					\label{itm:loosepath4}
					\item for each $i\in [\chi-1]$, the first and last element of $Q_i$ is flexible;
					\label{itm:loosepath5}
          \item
					every red loose path of length at most two connecting flexible elements from two different chains in $K$ contains a spine vertex;
					\label{itm:loosepath6}
					\item 
										\label{itm:loosepath7}
					for every subset $B\subseteq V(K)$ of size $\sigma$  and every collection of subsets of flexible vertices $A_1,\dots, A_{\chi-1}$ with $A_i\subseteq V(Q_i)\setminus B$ and $\size{A_i}\geq m$ for all $i\in[\chi-1]$, there exists a red edge $e\subseteq A_\ell\cup B$ for some $\ell\in [\chi-1]$ satisfying $1\leq \size{e\cap B}\leq k-1$.
				\end{enumerate}
We assume further that $Q_1,\dots, Q_{\chi-1}$ are chosen so that $\sum_{i \in [\chi-1]}v(Q_i)$ is maximal, that is,  $\sum_{i \in [\chi-1]} \ell_i$ is minimal.

Note that $\ell_i>0$ for all $i\in [\chi-1]$, as otherwise \cref{lem:chain_to_path} yields a red copy of $P_{n,1}^{(k)}$, contradicting our initial assumption.
Since $\ell_i \equiv 0 \pmod{k-1}$, it follows that $\ell_i\geq k-1$ for all $i\in [\chi-1]$. Let $W = V(K)\setminus \bigcup_{i\in [\chi-1]}V(Q_i)$. We then have
\begin{align}
    \size{W} &= \sum_{i \in [\chi-1]} (\ell_i-1) +c \label{eq:Wsize}\\
		&\geq (\chi-1)(k-2) +c . \label{eq:Wsize2}
\end{align}
By~\ref{itm:loosepath2} and \ref{itm:loosepath3}, each chain~$Q_i$ contains a flexible element~$S_i$ with at least 
\begin{align*}
\frac{2n/3 + 10k \ell_i}{1.2n/M} \ge M/2
\end{align*} flexible vertices; let $S^-_i$ denote the set of flexible vertices in $S_i$. 
Then $\size{S^-_i}\geq M/2$.

\begin{clm}\label{clm:red_matching}
    There is no red matching $\mathcal{M}$ in $K[W\cup \bigcup_{i \in [\chi-1]}S_i^-]$ such that $\size{V(\mathcal{M})\cap W} > (\chi-1)(k-2)$ and each edge $e\in \mathcal{M}$ satisfies $1\leq \size{e\cap W}\leq k-2$.
\end{clm}
\begin{proof}[Proof of \cref{clm:red_matching}]
    Suppose for a contradiction that such a red matching $\mathcal{M}$ exists. 
		Each edge $e\in \mathcal{M}$ contains at least two vertices outside of~$W$, so \ref{itm:loosepath6} implies that $e\setminus W\subseteq V(S_{i_e}^-)$ for some~$i_e\in [\chi-1]$. 
		By averaging, there exists an $i\in [\chi-1]$ and a submatching $\mathcal{M}'\subseteq \mathcal{M}$ such that $i_e = i$ for all $e\in \mathcal{M}'$ and $ \size{V(\mathcal{M}')\cap W} > k-2$.
		Moreover, by discarding extra edges, we may assume $ \size{V(\mathcal{M}')\cap W} \leq 2k-2$. 

    We now show that we can add the edges of~$\mathcal{M}'$ into the clique chain~$Q_i$. 
		Write $\mathcal{M}' = \set{e_1,\dots, e_p}$.
		Let $v_0$ be a spine vertex in~$S_i$ and $v_0'$ be the other one, if it exists, so that $S^-_i = S_i \setminus \{ v_0, v_0' \}$.
		Since each edge~$e_j \in \mathcal{M}'$ contains at least two vertices of~$S_i^-$ and $S_i^-$ induces a red clique of order at least $20k^2$ in~$K$, there are red edges $f_1,\dots, f_p$ in $K[S_i\setminus v_0']$ such that $v_0\in f_1$ and the edges $f_1,e_1,\dots,f_p,e_p$ form a red loose path~$P$ of length~$2p$ (in this order). 
		We may assume that~$P$ starts at $v_0$ and ends at~$v_1$.  
		Let $S_i'$ be a subset of $S_i\setminus (V(P) \setminus \set{v_1})$ of largest size containing both $v_1$ and $v_0'$ and satisfying $\size{S_i'}\equiv 1 \pmod{k-1}$. 

		Construct the open $1$-clique chain~$Q_i'$ from~$Q_i$ by removing the element~$S_i$ and inserting the elements~$f_1,e_1,\dots, f_p,e_p,S_i'$ in its place (possibly in reversed order if $S_i$ is the first element of~$Q_i$).
		Clearly, $Q_i'$ is red. 
		Observe that, if $S_i$ is the first (last) element of~$Q_i$, after this operation~$S_i'$ will be the first (last) element of~$Q_i'$.

		It remains to show that replacing~$Q_i$ with~$Q_i'$ gives a collection of red 1-clique chains also satisfying \ref{itm:loosepath1}--\ref{itm:loosepath7} but containing more vertices, which will contradict the choice of our initial collection. 
		Clearly \ref{itm:loosepath5} holds by construction, and~\ref{itm:loosepath6} and~\ref{itm:loosepath7} are inherited from $Q_1,\dots, Q_{\chi-1}$. 
    Note that 
    \begin{align*}
        v(Q_i') \geq v(Q_i) + \size{V(\mathcal{M}')\cap W} - (k-2) > v(Q_i).
    \end{align*}
    Hence, $v(Q_i') = n-\ell'_i$ for some $0 <\ell'_i < \ell_i - (k-1)$ and so \ref{itm:loosepath1} holds. 
				The number of flexible elements in~$Q_i'$ is the same as that in~$Q_i$ implying~\ref{itm:loosepath2}. 
		The number of flexible vertices is reduced by at most
    \begin{align*}
        2k\size{\mathcal{M}'} + (k-2) \leq 4k^2 \leq 10k(\ell_i - \ell_i'),
    \end{align*}
    so $Q_i'$ has at least $2/3n+10k\ell_i'$ flexible vertices, implying~\ref{itm:loosepath3}. 
		Finally, the flexible element $S_i'$ still contains at least $M/2-4k^2\geq M/4$ vertices implying~\ref{itm:loosepath4}.
\end{proof}

Consider a largest red matching~$\mathcal{M}$ in~$K[W\cup \bigcup_{i \in [\chi-1]}S_i^-]$ such that $1\leq \size{e\cap W}\leq k-2$ for each edge~$e\in \mathcal{M}$. 
Then, by \cref{clm:red_matching}, we have $\size{V(\mathcal{M})\cap W} \leq (\chi-1)(k-2)$.
Let $W' = W\setminus V(\mathcal{M})$ and $S'_i = S_i\setminus V(\mathcal{M})$ for all $i\in [\chi-1]$.
Then
\begin{align}
    \size{W'} &\geq \size{W} - (\chi-1)(k-2) \label{eq:W'size}\\
		& \overset{\mathclap{\text{\eqref{eq:Wsize2}}}}{\geq} c \ge \sigma, \label{eq:W'size2}\\
    \size{S_i'}&\geq \size{S_i} - \chi k^2 \geq M/4 \geq m +  \binom{\tau(k-1,\sigma)+1}{k-1} \label{eqn:S_i'size}
\end{align}
for all $i\in [\chi-1]$.
By the maximality of~$\mathcal{M}$, there is no red edge in $K[W'\cup \bigcup_{i \in [\chi-1]}S_i']$ with $1\leq \size{e\cap W'}\leq k-2$. 

\medskip
We now show that $\sigma \ge k-1$. 
Suppose to the contrary that $\sigma\leq k-2$.
Let $B$ be a subset of~$W'$ of size~$\sigma$, which exists by~\eqref{eq:W'size2}.
By~\ref{itm:loosepath7}, there exists a red edge in $K[B\cup \bigcup_{i \in [\chi-1]}S_i']$ satisfying $1\leq \size{e\cap W'}\leq \sigma \leq k-2$, a contradiction. Next, we consider two cases depending on the size of $W'$.

Suppose first that $\size{W'} > \tau(k-1,\sigma)$.
We may assume, by removing additional vertices if necessary, that $\size{W'} = \tau(k-1,\sigma)+1$. 
Recall that every red edge in $K[W'\cup \bigcup_{i \in [\chi-1]}S_i']$ satisfies $\size{e\cap W'}\in \set{0,k-1}$.
We now define a $(k-1)$-graph~$J$ on vertex set~$W'$ as follows.
We consider each $(k-1)$-subset $f\subseteq W'$ in turn and add it as an edge to~$J$ if there is a previously unused vertex $w_f\in \bigcup_{i\in [\chi-1]}S_i'$ such that $f\cup \set{w_f}$ is a red edge in~$K$; in this case, we mark $w_f$ as used.
By~\eqref{eqn:S_i'size} and~\ref{itm:loosepath7}, $\alpha(J) < \sigma$. 
By \cref{def:tau} and since  $v(J) = \size{W'} > \tau(k-1,\sigma)$, we know that $J$ contains a loose path of length two with edges~$f,f'$. 
This in turn leads to a red loose path of length two with edges $f\cup \set{w_f}$ and~$f'\cup \set{w_{f'}}$ in~$K$. 
By~\ref{itm:loosepath6},  $w_f$ and~$w_{f'}$ must belong to the same set~$S'_i$, and so we can extend the chain~$Q_i$ using~$f\cup \set{w_f},f'\cup \set{w_{f'}}$ in a similar fashion as in the proof of \cref{clm:red_matching}.
\medskip

Therefore we may assume that $\size{W'} \le \tau(k-1,\sigma)$.
Together with \eqref{eq:Wsize} and~\eqref{eq:W'size}, we have 
\begin{align*}
	\tau(k-1,\sigma) & \ge \size{W'} 
	\ge  \size{W} - (\chi-1)(k-2)  =  \sum_{i \in [\chi-1]} (\ell_i-1) + c - (\chi-1)(k-2) \\ 
	& \ge  \sum_{i \in [\chi-1]} (\ell_i-1) + \tau(k-1,\sigma) - 2k+3 - (\chi-1)(k-2) \\
	&= \sum_{i \in [\chi-1]} \ell_i + \tau(k-1,\sigma) - (\chi+1)(k-1)+1.
\end{align*}
Thus, we have $ \sum_{i \in [\chi-1]} \ell_i  \leq   (\chi+1)(k-1)-1$.
Since $\ell_i \equiv 0 \pmod{k-1}$ and $\ell_i >0$ for each $i\in [\chi-1]$, all but at most one $\ell_i$ are equal to $k-1$ and the remaining $\ell_i$ is at most~$2(k-1)$.
In this case, we let $S_i$ be the first element of the clique chain~$Q_i$ for each $i\in [\chi-1]$.
By~\ref{itm:loosepath4} and~\ref{itm:loosepath5}, the set~$S^-_i$ of flexible elements in~$S_i$ has size  at least $M/ 4-1$.
(We will no longer require that \ref{itm:loosepath4} holds, which is fine as we will only need to extend the $Q_i$ at most twice.)
Note that Claim~\ref{clm:red_matching} still holds. 
By~\eqref{eq:W'size2} and \ref{itm:loosepath7}, there exists some $j\in [\chi-1]$ and a red edge~$e$ in~$K[W'\cup \bigcup_{i \in [\chi-1]} S_i^-]$ such that $ e \cap S_j^-, e \cap W' \ne \emptyset$.
By Claim~\ref{clm:red_matching}, $\size{W'\cap e} = k-1$.
We add~$e$ to the chain~$Q_j$ as the first element. 
If $\ell_j = k-1$, then this together with~\cref{lem:chain_to_path} yields a contradiction. 
Otherwise, we can repeat this procedure with $S_j$ being the last element of~$Q_j$ instead. 
\end{proof}

We now sketch the proof of \cref{thm:loosecycles_upper}.
One can prove an analogue of Corollary~\ref{clm:chainspartition}, where each $Q_i$ is now closed, so \ref{itm:clm:chainspartition:firstlast} no longer applies, and \ref{itm:clm:chainspartition:crossing} is replaced by the following:
\begin{enumerate}[label={\rm (\alph*$'$)}, start=6]
	\item every pair of vertex-disjoint red loose paths of length at most two connecting flexible elements $S$ and $S'$ from two different chains in $K$ contains a spine vertex.\label{loose_cycle_f}
\end{enumerate}

Let $\mathcal{L}$ be a set of vertex-disjoint red paths of length at most two connecting two different chains which does not contain any spine vertices. By~\ref{loose_cycle_f}, we know that, for any pair of flexible elements $S$ and $S'$ from different chains, $\cL$ can contain at most one path connecting $S$ and $S'$. We now claim that, for any distinct $i,j\in [\chi-1]$, there can be at most two paths connecting~$V(Q_i)$ and $V(Q_j)$ in $\cL$. Indeed, assume there are at least three paths in $\cL$ connecting $V(Q_i)$ and $V(Q_j)$, and say that they connect the flexible elements $S_1,S_2$, and~$S_3$ in $Q_i$ and $S_1'$, $S_2'$, and $S_3'$ in~$Q_j$, respectively. Now, considering the six open subchains of $Q_i$ starting with some element $S_\alpha$ and ending with another element $S_\beta$, we conclude that at least one of those subchains contains at least $2v(Q_i)/3$ vertices; call this open subchain $Q_i'$. Without loss of generality, say $\alpha=1$ and $\beta=2$. Now, in a similar way, the corresponding flexible elements $S_1'$ and $S_2'$ in $Q_j$ bound a subchain $Q_j'$ containing at least half of the vertices of $Q_j$. We can then use the paths between~$S_1$ and $S_1'$ and $S_2$ and $S_2'$ to join the chains $Q_i'$ and $Q_j'$ into a closed clique chain containing at least $2v(Q_i)/3 +v(Q_j)/2 \geq (2/3 + 1/2)99n/100 > 1.1n$ vertices, giving a loose cycle of length~$n$. So we may assume that $\cL$ contains at most two paths connecting any pair of distinct chains.   Hence $\size{\cL} \leq 2\binom{\chi-1}{2}$. 
By setting $U = V (\mathcal{L} )$, we obtain the following the statement instead of~\ref{loose_cycle_f}:
\begin{enumerate}[label={\rm (\alph*$''$)},start=6]
	\item there exists a subset~$U\subseteq V(K)$ such that $\size{U} \le 4k\binom{\chi-1}{2}$ and every red loose path of length at most two connecting flexible elements from two different chains in $K \setminus U$ contains a spine vertex.
\end{enumerate}
The rest of the proof of~\cref{thm:loosecycles_upper} follows similarly the proof of~\cref{thm:loosepaths_upper}, where we remove the set~$U$ after defining $Q_1, \dots, Q_{\chi-1}$.

%%%%%%%%%%%%%%%%%%%%%%%%%%%%%%%%%%%%%%%%%%%
\section{Upper bound for $3$-uniform tight paths}\label{sec:transitive_ub}

In this section all hypergraphs are 3-uniform.
Recall that $TT_\chi$ is a transitive tournament on $[\chi]$ and $H(TT_\chi,m)$ is the tournament hypergraph associated to $TT_\chi$ in which each vertex class has exactly $m$ vertices (see \cref{def:tournament_hypergraph}). Recall that, given a 3-graph $F$ and subsets~${A,B, C\subseteq V(F)}$, we write~$E(A,B,C)$ for the set of  triples~$(a,b,c)\in A\times B\times C$ of distinct vertices such that~$abc\in E(F)$, and $d(A,B,C)$ for the ratio of $\size{E(A,B,C)}$ to the number triples of distinct vertices in $A\times B\times C$.

\subsection{Finding a large clique chain}

Our first lemma shows that, for any large collection of vertex-disjoint large red cliques in a red/blue-edge-coloured complete 3-graph, there will be enough blue edges to form a transitive tournament hypergraph or two of the red cliques will be connected by many (disjoint) red tight paths of length two. This result generalises Lemma 11 in~\cite{balogh2020ramsey}. 
 
\begin{lem}\label{lem:butterflies}
Let~$\chi,m\geq 2$ be integers and write~$R = \vec{R}(\chi)$. Then there exists an integer ${M_{\ref{lem:butterflies}} = M_{\ref{lem:butterflies}}(\chi,m)}$ such that the following holds for any~$M\geq M_{\ref{lem:butterflies}}$.
Let $K = K^{(3)}_{RM}$ be  red/blue-edge-coloured. Suppose $V(K)$ is partitioned into $R$ sets $V_1,\dots, V_R$, each inducing a red copy of~$K_M^{(3)}$. Then at least one of the following holds:
\begin{enumerate}[label=\normalfont(\alph*)]
    \item For any collection of subsets $W_i\subseteq V_i$ satisfying $\size{W_i}\geq M/2$ for all $i\in [R]$, there exists a red tight path of length two in $K[\bigcup_{i\in [R]}W_i]$ connecting $W_j$ and $W_{j'}$ for some distinct $j,j'\in [R]$. 
    \label{lem:butterflies:red}
    \item There exists a blue copy of~$H(TT_\chi,m)$.\label{lem:butterflies:blue}
\end{enumerate}
\end{lem}
\begin{proof}
Throughout the proof, we assume that~$M_{\ref{lem:butterflies}}$ is large and $M\geq M_{\ref{lem:butterflies}}$; we will specify how large~$M_{\ref{lem:butterflies}}$ needs to be in due course. We make no serious effort to optimise this constant.
Suppose~\ref{lem:butterflies:red} is false and consider a collection of subsets $W_i\subseteq V_i$ with $\size{W_i}\geq M/2$ for each~$i\in [R]$ such that no two distinct sets $W_i$ are connected by a red tight path of length two.
Thus, for each pair of vertices~$a\in W_i$ and~$b\in W_j$ with $i\neq j$, at least one of~$E(a,b, W_i)$ and~$E(a,b, W_j)$ is entirely blue.
Consider an auxiliary $R$-partite digraph $\vec{D}$ with vertex classes $W_1,\dots,W_R$ such that $(a,b)\in W_i\times W_j$ is an arc if and only if~$E(a,b, W_j)$ is entirely blue. 
 Let~$D$ denote the underlying (undirected) graph of~$\vec{D}$; note that~$D$ is a complete~$R$-partite graph with vertex classes~$W_1,\dots, W_R$.

We edge-colour $D$ so that an edge $ab$ with $a\in W_i$ and $b\in W_j$, where $i< j$, is coloured black if $(a,b)\in E(\vec{D})$ and white otherwise.
We now find a subset~$W'_i$ of~$W_i$ of size~$m$ so that each~$D[W'_i\cup W'_j]$ is monochromatic.
For~$t\in \mathbb{N}$, let~$B(t)$ denote the bipartite Ramsey number of~$K_{t,t}$, that is,~$B(t)$ is the minimum integer~$n$ such that in every red/blue-edge-colouring of~$K_{n,n}$ there exists a monochromatic copy of~$K_{t,t}$. It is well known that~$B(t)$ is finite for any~$t\geq 2$, see~\cite{beineke1976bipartite}\footnote{More precisely, we know that~$(1+o(1))(\sqrt{2}/e)t2^{t/2}\leq B(t) \leq (1+o(1))\log_2t2^{t+1}$ as~$t\to \infty$, see~\cite{conlon2008new,hattingh1998bipartite}.}.
We begin by setting~$W'_i = W_i$ for all~$i\in [R]$.
We now go through  all pairs of indices~$1\leq i < j\leq R$ and use the finiteness of bipartite Ramsey numbers to find large subsets of~$W'_i$ and~$W'_j$ such that $D[W'_i\cup W'_j]$ is monochromatic (and then discard all other vertices from the old~$W'_i$ and~$W'_j$). Provided that~$M/2\geq B(B(\dots B(B(m))))$, iterated~$\binom{R}{2}$ times, we can ensure that each~$W'_i$ has size at least~$m$ in the end. By discarding additional vertices, we may assume that each subset~$W'_i$ has size exactly~$m$. 

Define a tournament $T_R$ on $[R]$ such that $(i,j)$ with $i<j$ is an arc in $T_R$ if and only if $D[W'_i\cup W'_j]$ is black. Recall that a black edge $ab\in D[W'_i\cup W'_j]$ means that $(a,b)$ is an arc in~$\vec{D}$, implying that $E_K(a,b,W'_j)\subseteq E_K(a,b,V_j)$ is entirely blue (a similar statement holds if $ab$ is white). Hence $K$ contains a blue copy of $H(T_R,m)$ with vertex classes $W'_1,\dots,W'_R$. By the definition of $R$, it follows that $T_R$ contains a copy of $TT_{\chi}$, which in turn implies that $K$ contains a blue copy of $H(TT_{\chi}, m)$. 
\end{proof}

\begin{cor}\label{lem:long_tight_chain}
Let~$\chi,m\geq 2$ be integers, $\eps>0$, and~$R = \vec{R}(\chi)$. Then there exists an integer~$M_{\ref{lem:long_tight_chain}} = M_{\ref{lem:long_tight_chain}}(\chi,m, \eps)$ such that the following holds for any~$M\geq M_{\ref{lem:long_tight_chain}}$.
Let $t\geq 1$ be an integer and $K^{(3)}_{tM}$ be red/blue-edge-coloured with no blue copy of $H(TT_\chi, m)$. Suppose $V(K^{(3)}_{tM})$ is partitioned into $t$ sets $V_1,\dots, V_t$, each inducing a red copy of $K_M^{(3)}$.
Then there exists a red closed $2$-clique chain~$Q$ on at least $\frac{t(M-30)}{R-1}$ vertices such that each flexible element of $Q$ is of order at least $M/2$ and all but at most $\eps v(Q)$ vertices of $Q$ are flexible.
\end{cor}
\begin{proof}
    This follows directly by \cref{lem:alpha_chains,lem:butterflies} and taking the largest of the resulting chains.
\end{proof}

\subsection{Further auxiliary results}

After building a suitable clique chain, we will require several additional auxiliary results to complete the proof of \cref{thm:tight_paths_ub}. We begin with a simple proposition concerning the difference between consecutive values of~$\vec{R}(\chi)$. 

\begin{prop}\label{lem:consecutiveRp}
For any~$\chi\geq 3$, we have~$\vec{R}(\chi)\geq \vec{R}(\chi-1)+2$.
\end{prop}
\begin{proof}
    Let~$N = \vec{R}(\chi-1)+1$. We provide an orientation of~$K = K_N$ that contains no copy of~$TT_\chi$. Let~$v_1,v_2\in V(K)$ be arbitrary distinct vertices. By the definition of~$\vec{R}(\chi-1)$, we know that~$K-\set{v_1,v_2}$ has an orientation without a copy of~$TT_{\chi-1}$. We extend this orientation to all of~$K$ by adding the arcs~$(v_1,u)$ and~$(u,v_2)$ for every~$u\in V(K)\setminus\set{v_1,v_2}$; finally we add the arc~$(v_2,v_1)$.

    Suppose for a contradiction that there is a copy~$T$ of~$TT_\chi$. Since~$K-\set{v_1,v_2}$ was oriented without a copy of~$TT_{\chi-1}$, we know that~$T$ must contain both~$v_1$ and~$v_2$. However, for every $u\in V(K)\setminus\set{v_1,v_2}$, the vertices~$v_1,v_2,$ and~$u$ form a directed cycle, and thus cannot be part of a transitive tournament on at least three vertices. 
\end{proof}

Next, we show a simple random embedding lemma.

\begin{lem}\label{lem:random_embedding}
Let~$\chi,m\geq 2$ be integers. Then there exists a constant~$\gamma_{\ref{lem:random_embedding}} = \gamma_{\ref{lem:random_embedding}}(\chi,m)$ such that the following holds. Let~$F$ be a subhypergraph of a tournament hypergraph associated to $TT_{\chi}$ and~$V_1,\dots, V_{\chi}$ be its vertex classes. If~$\size{V_i}\geq 1/\gamma_{\ref{lem:random_embedding}}$ for each~$i\in [\chi]$ and~$d(V_i,V_i,V_j)\geq 1-\gamma_{\ref{lem:random_embedding}}$ for every arc~$(i,j)$ in  $TT_{\chi}$, then~$F$ contains a copy of~$H(TT_\chi,m)$.
\end{lem}
\begin{proof}

    Choose~$\gamma = \gamma_{\ref{lem:random_embedding}} = (\chi^2m^3)^{-1}$, and let~$F$ be a 3-graph satisfying the hypothesis. For each~$i\in [\chi]$, pick~$m$ vertices from~$V_i$ uniformly at random and independently from one another, allowing repetitions; let~$S_i$ be the (multi)set of chosen vertices. We will show that the probability that~$F[S_1\cup\dots\cup S_\chi]$ does not contain a copy of~$H(TT_{\chi},m)$ is less than 1, which then implies the desired result.

    Let~$(i,j)$ be an arc of~$TT_\chi$ and consider arbitrary vertices~$a,b\in S_i$ and~$c\in S_j$. The probability that~$a$ and~$b$ are the same is at most~$\frac{1}{\size{V_i}}\leq \gamma$, and the probability that~$a\neq b$ but~$abc$ is not an edge in~$F$ is at most~$\gamma$. The total number of triples~$abc$, where~$a,b\in S_i$ and~$c\in S_j$ and~$(i,j)$ is an arc of~$TT_{\chi}$ is bounded above by~$\binom{\chi}{2}m^3$. Taking a union bound over all triples, we find that the probability that any of them fails to be an edge in~$F$ is at most~$\binom{\chi}{2}m^3 2\gamma < \frac{\chi^2m^3}{2} \frac{2}{\chi^2m^3} = 1$. 
\end{proof}

A well-known result of Erd\H{o}s~\cite{erdos1964extremal} shows that any sufficiently large 3-graph with positive edge-density contains a copy of any fixed complete 3-partite 3-graph. A complete 3-partite 3-graph with the same number of vertices in each part  is easily seen to contain a spanning tight path. Our next lemma is in the same spirit; as we will use the precise dependencies among the constants, we include the short proof.

\begin{lem}\label{lem:absorbing_block}
    Let~$0<\eta\leq 1$ be a real number and ~$d\geq 1$ be an integer. Let~$F$ be a \mbox{$3$-graph} and~$A,B\subseteq V(F)$ be such that~$\size{A}\geq 4d\frac{e^d}{\eta^d}$,~$\size{B}\geq \frac{d}{\eta}$, and~$A\cap B = \emptyset$. Further, assume that~$A$ induces a clique in~$F$ and that~$d(A,A,B)\geq \eta$.  Then $F$ contains a tight path of the form~$a_1a_2b_1a_3a_4\dots a_{2d-1}a_{2d}b_ta_{2d+1}a_{2d+2}$ with $a_i\in A$ and $b_j\in B$ for all $i\in [2d+2]$ and $j\in [d]$.
\end{lem}

\begin{proof}
    We count the number of pairs~$(\set{a_1,a_2}, \set{b_1,\dots, b_d})$ such that~$a_1,a_2\in A$, $b_1,\dots,b_d\in B$ (with all~$d+2$ vertices being distinct), and~$a_1a_2b_i\in F$ for all~$i\in [d]$. 
    The number of pairs as described above equals
    \begin{align*}
        \frac12\sum\limits_{\substack{a_1,a_2\in A\\a_1\neq a_2}} \binom{\size{E(a_1, a_2, B)}}{d} &\geq \binom{|A|}{2} \binom{\sum\limits_{\substack{a_1,a_2\in A\\a_1\neq a_2}}\frac{\size{E(a_1, a_2, B)}}{\size{A}(\size{A}-1)}}{d}\\
        & = \binom{\size{A}}{2} \binom{d(A,A,B)\size{B}}{d} \geq \binom{\size{A}}{2} \binom{\eta\size{B}}{d},
    \end{align*}
    where the first step comes from Jensen's inequality and the last follows from our assumption on the density~$d(A,A,B)$.
    Therefore, by averaging, there exists a set~$\set{b_1,\dots, b_d}\subseteq B$ of size~$d$ such that the number of pairs~$(\set{a_1,a_2},\set{b_1,\dots, b_d})$ of the required form
    is at least 
    \begin{align}\label{eq:num_edges_G}
        \frac{\binom{\size{A}}{2} \binom{\eta\size{B}}{d}}{\binom{\size{B}}{d}} &\geq \binom{\size{A}}{2} \frac{\parens*{\frac{\eta\size{B}}{d}}^d}{\parens*{\frac{e\size{B}}{d}}^d} = \binom{\size{A}}{2}\frac{\eta^d}{e^d} > d\size{A},
    \end{align}
    where we used standard estimates on the binomial coefficient and the facts that~$\eta\size{B}\geq d$ and $\size{A}\geq 4d\frac{e^d}{\eta^d}$.

    Now, we build an auxiliary graph~$G$ on vertex set~$A$, in which~$a_1a_2$ is an edge whenever~$a_1a_2b_i\in F$  for all~$i\in [d]$. By~\eqref{eq:num_edges_G}, $G$ has more than $d\size{A}$ edges. Thus, by a classical result of Erd\H{o}s and Gallai~\cite{gallai1959maximal}, $G$ contains a path $a_1a_2\dots a_{2d+2}$ of length $2d+1$. Then~$a_ia_{i+1}b_j\in F$ for all~$i\in [2d+1]$ and~$j\in [d]$, implying that~$a_1a_2b_1a_3a_4\dots a_{2d-1}a_{2d}b_da_{2d+1}a_{2d+2}$ is indeed a tight path. 
\end{proof}

\subsection{Proof of~\cref{thm:tight_paths_ub}}

We are now ready to prove \cref{thm:tight_paths_ub}. 
We first handle the case~$\chi=2$, which requires a slightly different argument, and then we move on to the proof for~$\chi\geq 3$.

\begin{lem}\label{lem:p2case}
    For any~$\eps>0$ and $m\geq 2$, there exists an integer~$n_{\ref{lem:p2case}} = n_{\ref{lem:p2case}}(\eps,2,m)$ such that, for all~$n\geq n_{\ref{lem:p2case}}$, we have $
        R(P_{n,2}^{(3)},H(TT_2,m))\leq \parens*{1+\eps}n$.
\end{lem}
\begin{proof} 
Let~$M = \max\set*{M_{\ref{lem:long_tight_chain}}(2,m), \frac{200}{\eps}}$ and $n_{\ref{lem:p2case}} = \ceil{ 2\eps^{-1}R(K_M^{(3)}, K_{2m}^{(3)})}$.
Let $n\geq n_{\ref{lem:p2case}}$, $N = \parens*{1+\eps}n$, and $K = K_N^{(3)}$. Suppose there exists a~$(P^{(3)}_{n,2},H(TT_2,m))$-free edge-colouring of~$K$. By~\cref{lem:monochromatic_cliques}, we can find vertex sets~$V_1,\dots, V_t$ covering all but at most~$\eps n/2$ vertices of $K$, each inducing a red copy of~$K_M^{(3)}$ or a blue copy of~$K_{2m}^{(3)}$. As~$H(TT_2,m)\subseteq K_{2m}^{(3)}$, each~$V_i$ induces a red copy of~$K_M^{(3)}$.

Observe that~$t\geq \frac{N-\eps n/2}{M} =\parens*{1+\frac{\eps}{2}}\frac{n}{M}$. Then by \cref{lem:long_tight_chain}, there exists a red closed $2$-clique chain $Q$ on at least $t(M-30) = \parens*{1+\frac{\eps}{2}}n - \frac{30}{M}\parens*{1+\frac{\eps}{2}} n \geq n+1$ vertices. An application of \cref{lem:chain_to_path} then implies the existence of a red tight path on $n$ vertices. 
\end{proof}

Finally, we prove the general case of our main theorem, using the case~$\chi = 2$ as the base for an inductive proof. Recall that, given a red/blue-edge-coloured 3-graph $F$ and subsets ${A,B, C\subseteq V(F)}$, we write $d_r(A,B,C)$ and $d_b(A,B,C)$ for the  density of the red and blue $ABC$-edges, respectively.

\begin{proof}[Proof of~\cref{thm:tight_paths_ub}]
We proceed by induction on~$\chi$. \cref{lem:p2case} gives the base case~$\chi=2$.
Now suppose~$\chi\geq 3$ and the theorem holds for~$\chi-1$ and any~$m\geq 2$. Let~$\eps> 0$ and~$m\geq 2$ be arbitrary.  Before delving into the details of the proof, we give a brief overview. As in the proof of \cref{thm:loosepaths_upper}, we proceed by contradiction, assuming the existence of a~$(P^{(3)}_{n,2},H(TT_{\chi},m))$-free red/blue-edge-colouring of some $K_N^{(3)}$. We then apply~\cref{lem:long_tight_chain} to find \emph{one} red 2-clique chain~$Q$ containing slightly more than $2n/3$ vertices. Our goal then is to absorb more vertices into this chain. For this, we go through the flexible elements of $Q$ one by one and, for each flexible element $S$, we try to find a tight path containing most vertices of $S$ and about half as many vertices not previously contained in the chain. By using the induction hypothesis, we find a 
 blue copy of $H(TT_{\chi-1}, q)$ for some large constant $q$ that is vertex-disjoint from $Q$; between $S$ and this copy, we either find many blue edges, giving us a blue copy of $H(TT_{\chi}, m)$, or we find enough red edges to build the required tight path. Eventually, unless we find a blue copy of $H(TT_{\chi}, m)$, we succeed in building a tight path on $n$ vertices, contradicting our initial assumption.

Without loss of generality, we assume that $\eps < 1/10$. We begin by fixing some constants; let 
\begin{align*}
    \gamma &= \gamma_{\ref{lem:random_embedding}}(\chi,m),\quad d>\ceil{10/\eps}, \quad R = \vec{R}(\chi),\quad q = \ceil*{d/\gamma},\\
    M &= \max\set*{M_{\ref{lem:long_tight_chain}}(\chi,m), \ceil*{\frac{200k}{\eps}},\ceil*{\frac{16de^d}{\eps\gamma^d}},\ceil*{\frac{4}{\eps\gamma}}},\\
    n_0(\chi,m,\eps) &= \max \set*{n_0(\chi-1,q,\eps)+1, \ceil*{R(K_M^{(3)}, K_{\chi m}^{(3)})/\eps}}.
\end{align*}
Let $n\geq n_0$ and~$N =\parens*{\frac{2}{3}+\eps}(R-1)n$. Suppose there exists a~$(P^{(3)}_{n,2},H(TT_{\chi},m))$-free red/blue-edge-colouring of~$K=K_N^{(3)}$. Let~$V_1,\dots, V_t$ be the sets guaranteed by~\cref{lem:monochromatic_cliques}, covering all but at most~$\eps n$ vertices of~$K$ and each inducing a red copy of $K_M^{(3)}$ or a blue copy of $K_{\chi m}^{(3)}$; since~$H(TT_{\chi},m)\subseteq K_{\chi m}^{(3)}$, we know that each~$V_i$ must induce a red copy of~$K_M^{(3)}$. 

By~\cref{lem:long_tight_chain} applied to the graph $K[\bigcup_{t\in [t]}V_i]$ and our assumption that there is no blue copy of $H(TT_\chi,m)$, there exists a red closed $2$-clique chain $Q$ in which each flexible element contains at least $M/2$ vertices, all but $\eps v(Q)$ vertices are flexible, and
\begin{align}\label{eq:vQ}
    v(Q) \geq t\frac{M-30}{R-1} \geq \frac{(N-\eps n)(M-30)}{M(R-1)} 
    = \left( \frac{2}3 +\frac{R - 2}{R- 1}\eps \right) n  \cdot \frac{M-30}{M}    
    \geq
        (1+0.9\eps)\frac23n,
\end{align}
where the last inequality follows since $R\geq 4$, and hence $(R-2)/(R-1)\geq 2/3$, and $M\geq 200k/\eps$.
Observe that we may assume that~$v(Q) < n$; otherwise we are done by \cref{lem:chain_to_path}. 

Let $S_1,\dots, S_c$ be the flexible elements of $Q$, so that 
\begin{align}\label{eq:flexible_elements}
    \sum_{i\in [c]} \size{S_i} \geq (1-\eps)v(Q).
\end{align}
Our goal now is to absorb more vertices into this clique chain. For this, we will consider each flexible element~$S_\ell$ in turn and replace it by a red tight path containing more vertices.

Note that $N-n = \parens*{\frac23+\eps}\parens*{ \vec{R}(\chi) - \frac{3}{2+3\eps} - 1}n$. If $\chi = 3$, then the fact that~$\vec{R}(3) = 4$ implies that this quantity is at least $\parens*{1+\eps}n$. Otherwise, by \cref{lem:consecutiveRp}, we have~$\vec{R}(p)-2\geq \vec{R}(p-1)$, and so we conclude that 
\begin{align}\label{eq:induction}
    N-n \geq \parens*{\frac23+\eps} \parens*{\vec{R}(\chi-1)-1}n \geq R(P_{n,2}^{(3)}, H(TT_{\chi-1}, q)),    
\end{align}
by the induction hypothesis. 

\begin{clm}\label{clm:absorb}
    For every $\ell\in[c] $, there exists a $2$-clique chain $Q^\ell$ on $v(Q)+\sum\limits_{i=1}^{\ell}\frac{d}{2d+2}(1-\eps)\size{S_i}$ vertices whose flexible elements are precisely $S_{\ell+1},\dots, S_c$.
\end{clm}

\begin{proof}[Proof of \cref{clm:absorb}]
    We set $Q^0 = Q$.  Let $\ell\in [c]$ and suppose we have found the required clique chain $Q^{\ell-1}$. Consider the set $S_\ell$, which forms a flexible element in $Q^{\ell-1}$; let $x_1,x_2,y_1$, and $y_2$ be the spine vertices (with respect to $Q^{\ell-1}$) contained in $S_\ell$, where $x_1$ and $x_2$ are contained in the previous element of $Q^{\ell-1}$ and $y_1$ and $y_2$ are contained in the next element (if either of these elements is not defined, choose the corresponding vertices arbitrarily). 
    
    Set $I = S_\ell\setminus\set{x_1,x_2,y_1,y_2}$ and $W = (V(K)\setminus V(Q^{\ell-1}))\cup I$, and note that $\size{I}\geq M/4$.
    Consider a longest red tight path $P$ in $K[W]$ whose first two and last two vertices are all in~$I$ and 
    \begin{align*}
        d\size{V(P)\cap I} = (2d+2)\size{V(P)\setminus I}.
    \end{align*}
    We will show that $\size{V(P)\cap I}\geq (1-\eps)\size{S_\ell}$. 
    
    Suppose otherwise and let $A = I\setminus V(P)$ so that $\size{A} \geq \eps\size{S_\ell} \geq \eps M/4$. Since $K$ has no red tight path on $n$ vertices,  \cref{lem:chain_to_path} implies that $v(Q^{\ell-1}) < n$.
    Set $K^- = K-V(Q^{\ell-1})$ and note that $v(K^-) > N-n \geq R(P_{n,2}^{(3)}, H(TT_{\chi-1}, q))$ by \eqref{eq:induction}. Hence, there exists a blue copy of $H(TT_{\chi-1},q)$ in $K^-$; let $Y_1,\dots, Y_{\chi-1}$ denote its vertex classes. Recall that $q\geq \gamma^{-1}$ and $\size{A} \geq \eps M/4 \geq \gamma^{-1}$. Now, if $d_b(A,A,Y_j) \geq 1-\gamma$ for all $j\in [\chi-1]$, then by \cref{lem:random_embedding}, there is a blue copy of $H(TT_{\chi}, m)$ in $K$ with vertex set contained in $A\cup Y_1\cup\dots\cup Y_{\chi-1}$, a contradiction. 
    
    Therefore, there exists some~$j\in [\chi-1]$ such that~$d_r(A,A,Y_j) \geq \gamma$. Then, by~\cref{lem:absorbing_block}, since~$\size{A}\geq \eps M/4\geq 4de^d/\gamma^d$ and~$\size{Y_j} = q\geq d/\gamma$ we can find vertices~$a_1,\dots,a_{2d+2}\in A$ and \mbox{$b_1,\dots,b_d\in Y_j$} such that~$a_1a_2b_1a_3a_4b_2\dots a_{2d-1}a_{2d}b_da_{2d+1}a_{2d+2}$ is a red tight path.  Attaching this segment to the end of $P$ results in a longer tight path containing $2d+2 + \size{V(P)\cap I}$ vertices of $I$ and $d + \size{V(P)\setminus I}$ vertices not in $I$, thus preserving the required ratio of vertices in $I$ and not in $I$. Moreover, the first two and the last two vertices of this new path are all in $I$, contradicting the maximality of $P$. 

    Thus, $\size{V(P)\cap I}\geq (1-\eps)\size{S_\ell}$, as claimed, and therefore $\size{V(P)\setminus I} \geq  \frac{d}{2d+2}(1-\eps)\size{S_\ell}$. Let~$P'$ be a tight path obtained from $P$ by adding the vertices in $I\setminus V(P)$ to the end of $P$.  
    We now construct the clique chain $Q^{\ell}$ by replacing the flexible element $S_\ell$ in~$Q^{\ell-1}$ by~$x_1x_2 P' y_1y_2$. 
\end{proof}

\medskip
Let $Q^c$ be the final path guaranteed by \cref{clm:absorb}; then 
\begin{align*}
    v(Q^c)&\geq v(Q)+\sum\limits_{i\in [c]}\frac{d}{2d+2}(1-\eps)\size{S_i} \geq v(Q)+\frac{d}{2d+2}(1-\eps)^2v(Q)  \\
    & \geq (1+0.9\eps)\frac{2}{3}n \parens*{1+\frac{d}{2d+2}(1-\eps)^2} \geq n,
\end{align*}
where the second and last inequalities come from  \eqref{eq:flexible_elements} and  \eqref{eq:vQ}, respectively.
Hence, by \cref{lem:chain_to_path}, there exists a red tight path on~$n$ vertices, contradicting our initial assumption.
\end{proof}

%%%%%%%%%%%%%%%%%%%%%%%%%%%%%%%%%%%%%%%%%%%

\section{Conclusion and open problems}\label{sec:conclusion}
In this paper, we attempted to extend Burr's result showing that every sufficiently long path is $H$-good  for every graph $H$ (\cref{thm:burr}) to higher uniformities. In the hypergraph setting, we  focused on $\ell$-paths. Our work leads to several natural possible future directions.

First, our methods for showing upper bounds work well when $\ell\in \set{1,k-1}$, but cannot be applied when $2\leq \ell\leq k-2$ without some new ideas. The issue arises in the construction of our clique chains: in the case $\ell=k-1$, a tight path of length $q_0$ connecting two disjoint sets does not ``spill out'' of those sets; on the other hand, when $\ell=1$, we also have $q_0=1$, and an edge connecting two disjoint sets that intersects a third set can be viewed as connecting every pair of the three sets. For other values of $\ell$, we do not have as much control over where the connecting paths lie.

While our bounds for loose paths are best possible, this is not the case for loose cycles (see Theorem~\ref{thm:loosecycles_upper} and Proposition~\ref{prop:loose_cycle_lower}). It would be interesting to determine the correct additive constant. 

In the case of tight paths, it would be interesting to extend our work, in particular our positive result, to uniformities larger than three. Our methods do not easily extend to larger $k$, and it is also not clear how to generalise the class of transitive tournament hypergraphs. Another natural problem is to turn \cref{thm:tight_paths_ub} into an exact result. Additionally, our examples proving \cref{prop:non_transitive} are fairly dense, leading us to ask: how sparse can we make those examples and how large can the ratio $R(P_{n,2}^{(3)}, H)/n$ get? Note that, given a fixed $k$-graph $H$, the disjoint union of $H$ and a tight path (or more generally, any $\ell$-path) on $n$ vertices is a bounded-degree hypergraph. It was shown by Cooley, Fountoulakis,  K{\"u}hn,  and Osthus~\cite{cooley2009embeddings}, following earlier results for 3-graphs that, for any $k\geq 3$ and $\Delta\geq 1$, there exists a constant $C = C(k,\Delta)$ such that any $m$-vertex $k$-graph $F$ with maximum (vertex) degree bounded above by  $\Delta$ satisfies $R(F,F) \leq Cm$. 

Finally, it would be interesting to investigate whether there is a more suitable notion of Ramsey goodness for hypergraphs, e.g.,
to find another natural lower bound on the Ramsey number of a pair of hypergraphs. 

\section*{Acknowledgements}
We are grateful to Shagnik Das for the argument showing the upper bound on $\tau(k,\alpha)$ in \cref{prop:tau_bounds}. We also thank the anonymous referee for their helpful feedback.

\bibliographystyle{amsplain}
\bibliography{biblio}

\end{document}